\pdfoutput=1
\RequirePackage{ifpdf}
\ifpdf 
\documentclass[pdftex]{sigma}
\else
\documentclass{sigma}
\fi

\numberwithin{equation}{section}

\newtheorem{Theorem}{Theorem}[section]
\newtheorem{Corollary}[Theorem]{Corollary}
\newtheorem{Lemma}[Theorem]{Lemma}
\newtheorem{Proposition}[Theorem]{Proposition}
 { \theoremstyle{definition}
\newtheorem{Definition}[Theorem]{Definition}
\newtheorem{Example}[Theorem]{Example}
\newtheorem{Remark}[Theorem]{Remark} }

\newcommand{\knabla}{%
 \setbox0=\hbox{$\triangle $}%
 \raisebox{\ht0}{\rotatebox{180}{$\triangle $}}}

\begin{document}
\allowdisplaybreaks

\newcommand{\arXivNumber}{1511.01608}

\renewcommand{\thefootnote}{}

\renewcommand{\PaperNumber}{110}

\FirstPageHeading

\ShortArticleName{Flat Structure on the Space of Isomonodromic Deformations}

\ArticleName{Flat Structure on the Space\\ of Isomonodromic Deformations\footnote{This paper is a~contribution to the Special Issue on Primitive Forms and Related Topics in honor of~Kyoji Saito for his 77th birthday. The full collection is available at \href{https://www.emis.de/journals/SIGMA/Saito.html}{https://www.emis.de/journals/SIGMA/Saito.html}}}

\Author{Mitsuo KATO~$^\dag$, Toshiyuki MANO~$^\ddag$ and Jiro SEKIGUCHI~$^\S$}

\AuthorNameForHeading{M.~Kato, T.~Mano and J.~Sekiguchi}

\Address{$^\dag$~Department of Mathematics, College of Educations, University of the Ryukyus, Japan}
\EmailD{\href{mailto:mkato@nirai.ne.jp}{mkato@nirai.ne.jp}}

\Address{$^\ddag$~Department of Mathematical Sciences, Faculty of Science, University of the Ryukyus, Japan}
\EmailD{\href{mailto:tmano@math.u-ryukyu.ac.jp}{tmano@math.u-ryukyu.ac.jp}}

\Address{$^\S$~Department of Mathematics, Faculty of Engineering,\\
\hphantom{$^\S$}~Tokyo University of Agriculture and Technology, Japan}
\EmailD{\href{mailto:sekiguti@cc.tuat.ac.jp}{sekiguti@cc.tuat.ac.jp}}

\ArticleDates{Received March 19, 2020, in final form October 21, 2020; Published online November 03, 2020}

\Abstract{Flat structure was introduced by K.~Saito and his collaborators at the end of 1970's. Independently the WDVV equation arose from the 2D topological field theory. B.~Dubrovin unified these two notions as Frobenius manifold structure. In this paper, we study isomonodromic deformations of an Okubo system, which is a special kind of systems of linear differential equations. We show that the space of independent variables of such isomonodromic deformations can be equipped with a~Saito structure (without a metric), which was introduced by C.~Sabbah as a generalization of Frobenius manifold. As its consequence, we introduce flat basic invariants of well-generated finite complex reflection groups and give explicit descriptions of Saito structures (without metrics) obtained from algebraic solutions to the sixth Painlev\'e equation.}


\Classification{34M56; 33E17; 35N10; 32S25}

\renewcommand{\thefootnote}{\arabic{footnote}}
\setcounter{footnote}{0}

\vspace{-2mm}

\section{Introduction}

At the end of 1970's, K.~Saito introduced the notion of flat structure in order to study the structure of universal unfolding of isolated hypersurface singularities.
Independently the WDVV equation (Witten--Dijkgraaf--Verlinde--Verlinde equation)
arose from the 2D topological field theory \cite{DVV,Wi}.
B.~Dubrovin unified both the flat structure formulated by K.~Saito and the WDVV equation
 as Frobenius manifold structure.
Dubrovin not only formulated the notion of Frobenius manifold but also studied its relationship
with isomonodromic deformations of linear differential equations with certain symmetries.
Particularly
he derived a one-parameter family of Painlev\'e VI equation
from three-dimensional massive (i.e., regular semisimple) Frobenius manifolds.
Since then,
there are several generalizations of Frobenius manifolds such as $F$-manifold by C.~Hertling and Yu.~Manin \cite{He, HM}
and Saito structure (without metric) by C.~Sabbah~\cite{Sab}.
Concerning the relationship with the Painlev\'e equation,
A.~Arsie and P.~Lorenzoni \cite{AL0,Lo} showed that three-dimensional regular semisimple bi-flat $F$-manifolds are parameterized
by solutions to the (full-parameter) Painlev\'e VI equation,
which is regarded as an extension of Dubrovin's result.
Furthermore Arsie--Lorenzoni \cite{AL1} showed that three-dimensional regular non-semisimple bi-flat $F$-manifolds are parameterized
by solutions to the Painlev\'e V and IV equations.

The theory of linear differential equations on a complex domain is a classical branch of mathematics.
In recent years there has been a great progress in this branch.
One of its turning points is
the introduction of the notions of {\it middle convolution} and {\it rigidity index} by N.M.~Katz~\cite{Katz}.
With the help of his idea,
T.~Oshima \cite{OshiB, Oshi1,Oshi2} developed a classification theory of Fuchsian differential equations
in terms of their rigidity indices and spectral types;
\cite{HaBook} is a nice introductory text on the ``Katz--Oshima theory''.
(From the viewpoint of the theory of isomonodromic deformations,
the number of accessory parameters of an irreducible linear differential equation is equal to two plus the index of rigidity.)
{\it Okubo system} (which was systematically studied by K.~Okubo in~\cite{Ok})
plays a central role in these developments (cf.\ \cite{DR1,DR2,Oshi2,Yo}):
a matrix system of linear differential equations with the form
 \begin{equation} \label{eq:introOkubo}
 (z-T)\frac{{\rm d}Y}{{\rm d}z}=BY,
 \end{equation}
 where $T$ and $B$ are constant square matrices,
 is an Okubo system if $T$ is diagonalizable.
 We remark that P.~Boalch~\cite{Bo0} introduced and used the operation of convolution which is equivalent to Katz's
 middle convolution,
where it was derived in a very simple fashion from a scalar shift considered explicitly by Balser--Jurkat--Lutz \cite{BJL}.
(Some of viewpoints in the present paper are already found in the work \cite{BoC,Bo0,Bosix} of Boalch.
See Corollary~\ref{cor:corgexWDVVandPVI} and Section~\ref{sec:algsol}.)

In this paper, we study isomonodromic deformations of an Okubo system (\ref{eq:introOkubo}).
We show that the space of independent variables of isomonodromic deformations of an Okubo system
can be equipped with a Saito structure (without a metric).
(In the sequel, we abbreviate ``Saito structure (without metric)'' to ``Saito structure'' for brevity.)
As its consequence, we obtain the following results:
\begin{itemize}\itemsep=0pt
 \item[(I)] introduction of flat coordinates on the orbit spaces of well-generated finite complex reflection groups
 (Section~\ref{sec:complexreflection}),
 \item[(II)] explicit descriptions of Saito structures corresponding to algebraic solutions to the Painlev\'e VI equation
 (Section~\ref{sec:algsol}).
\end{itemize}
K.~Saito \looseness=-1 and his collaborators \cite{Sai2,SYS} defined and constructed flat coordinates on the orbit spaces of finite real reflection groups.
To extend them to finite complex reflection groups has been a~long-standing problem.
(I)~gives an answer to this problem for well-generated complex reflection groups.
Recently the construction of flat coordinates on the orbit spaces of well-generated complex reflection groups is crucially applied
in \cite{HMRS} to prove freeness of multi-reflection arrangements of such groups.
(II)~provides many new concrete examples of three-dimensional algebraic Saito structures that are not Frobenius manifolds,
which would be the first step toward classification of three-dimensional algebraic Saito structures and/or algebraic $F$-manifolds (cf.~\cite{He}).
(In this paper, the term ``flat structure'' stands for the same meaning as ``Saito structure''
because it is a natural extension of K. Saito's flat structure,
which would be justified by the result~(I) above.)

This paper is constructed as follows.
In Section~\ref{sec:okuboiso}, we introduce the notion of {\it extended Okubo system}
as a completely integrable Pfaffian system extending an Okubo system.
We see that an extended Okubo system is equivalent to an isomonodromic deformation of an Okubo system.
We observe that an extended Okubo system naturally induces a {\it Saito bundle} \cite{DH, Sab} on a complex domain,
which will be used in Section~\ref{subsec:Okubo}.

In Section~\ref{sec:logvecf}, we study the structure of logarithmic vector fields along a~divisor defined by a~monic polynomial of degree~$n$
\begin{equation} \label{eq:h(x)poly}
 h(x)=h(x',x_n)=x_n^n-s_1(x')x_n^{n-1}+\cdots+(-1)^ns_n(x'),
\end{equation}
where $x'=(x_1,\dots,x_{n-1})$ and each $s_i(x')$ is holomorphic with respect to $x'$,
which appears as the defining equation of the singular locus of an extended Okubo system.
In K.~Saito's construction of flat structures on the orbit spaces of real reflection groups,
the fact that the discriminants of real reflection groups have the form of~(\ref{eq:h(x)poly}) was used crucially (cf.~\cite{Sai2}).
The results in this section are used in Section~\ref{sec:complexreflection} in order to generalize Saito's construction.
In Section~\ref{subsec:Saito}, we review the general theory of Saito structure.
Particularly we introduce an extension of the WDVV equation,
which is not mentioned in~\cite{Sab} explicitly.
A solution to the extended WDVV equation is called a potential vector field,
which (at least locally) completely describes a Saito structure.
If a Saito structure admits a Frobenius metric, the potential vector field is integrated to a prepotential,
which is a~solution to the WDVV equation.
In Section~\ref{subsec:Okubo},
we show that the space of variables of an extended Okubo system can be equipped with
a Saito structure under some generic condition (Theorem~\ref{saitojacobian}).
Theorem~\ref{saitojacobian} is well known (at least) in the case of Frobenius manifolds.
The arguments here closely follow \cite[Chapter~VII]{Sab}:
we find a condition for that the Saito bundle induced by an extended Okubo system in Section~\ref{sec:okuboiso}
has a {\it primitive section}.
As a consequence of Theorem~\ref{saitojacobian}, we find that there is a correspondence between generic solutions satisfying some semisimplicity condition
 to the three-dimensional extended WDVV equation
and generic solutions to the Painlev\'e~VI equation (Corollary~\ref{cor:corgexWDVVandPVI}).
In the proof of Corollary~\ref{cor:corgexWDVVandPVI},
we use the fact that isomonodromic deformations of a rank-three Okubo system are governed by generic solutions to the Painlev\'e~VI equation,
which was treated by P.~Boalch~\cite{BoC,Bo0,Bosix} in order to construct algebraic solutions to the Painlev\'e VI equation
related with rank-three finite complex reflection groups.
(The operation to construct a $3\times 3$ Okubo system from a $2\times 2$ system by Boalch is equivalent to
Katz's middle convolution, but given in a~different approach from Katz~\cite{Katz} or Dettweiler--Reiter~\cite{DR1}.)
In Section~\ref{sec:complexreflection}, we treat the problem on the existence of a flat coordinate system
on the orbit space of an irreducible complex reflection group.
In the case of a real reflection group, K.~Saito~\cite{Sai2} (see also~\cite{SYS}) proved the existence of a flat coordinate system
on the orbit space based on the existence of a flat invariant metric on the standard representation space.
In this paper, instead of the flat invariant metric, we use a special kind of extended Okubo systems
called {\it $G$-quotient system}, whose fundamental system of solutions consists of derivatives by logarithmic vector fields
of linear coordinates on the standard representation space
and its monodromy group is isomorphic to the finite complex reflection group $G$ (see Remark~\ref{rem:gs vs okubo}).
Applying Theorem~\ref{saitojacobian} to the $G$-quotient system,
we introduce flat basic $G$-invariants (Theorem~\ref{prop:C_1}).
As its consequence, we see that the potential vector fields corresponding to well-generated complex reflection groups
have polynomial entries.
It is underlined that our proof of Theorem~\ref{prop:C_1} is constructive:
Theorem~\ref{prop:C_1} contains a method
to construct the flat basic invariants and the potential vector field explicitly
for each well-generated complex reflection group.
We explain it for an example (Example~\ref{ex:G26}).
In Section~\ref{sec:algsol}, we treat algebraic solutions to the Painlev\'e VI equation
as an application of Corollary~\ref{cor:corgexWDVVandPVI}.
Algebraic solutions to the Painlev\'e~VI equation were studied and constructed by
 many authors including N.J.~Hitchin \cite{Hit1,Hit2},
 B.~Dubrovin \cite{Du}, B.~Dubrovin and M.~Mazzocco~\cite{DM},
P.~Boalch \cite{BoC,Bo0,Bo1,Bo3,Bo2},
A.V.~Kitaev \cite{Kit1,Kit2}, A.V. Kitaev and R.~Vid\=unas \cite{KV,VK}, K.~Iwasaki~\cite{Iwa}.
In his review article~\cite{Bo4}, Boalch classified all the known algebraic solutions to the Painlev\'e~VI equation.
After that, Lisovyy and Tykhyy~\cite{LT} proved that Boalch's list was complete, that is,
they confirmed that there is no other solution than those in Boalch's list.
We give some examples of potential vector fields corresponding to algebraic solutions to the Painlev\'e~VI equation.
Related results are given in~\cite{KMS4}.

\section{Okubo system and its isomonodromic deformation} \label{sec:okuboiso}

Let $T$ and $B_{\infty}$ be $(n\times n)$-matrices.
If $T$ is a diagonalizable matrix,
the system of ordinary linear differential equations
\begin{equation} \label{eq:okubo}
 (zI_n-T)\frac{{\rm d}Y}{{\rm d}z}=-B_{\infty}Y
\end{equation}
is called an {\it Okubo system}~\cite{Ok}.
We extend (\ref{eq:okubo}) to a completely integrable Pfaffian system of the form
\begin{equation} \label{eq:okubopfaff}
{\rm d}Y=-(zI_n-T)^{-1}\bigl({\rm d}z+\tilde{\Omega}\bigr)B_{\infty}Y
\end{equation}
on $\mathbb{P}^1\times U$, where $U$ is a simply-connected domain in~$\mathbb{C}^n$
and $\tilde{\Omega}$ is an $(n\times n)$-matrix valued $1$-form satisfying the conditions below.
We put the following assumptions on~(\ref{eq:okubopfaff}):
\begin{enumerate}\itemsep=0pt
 \item[(A1)] $T$ (resp.\ $\tilde{\Omega}$) is a diagonalizable $(n\times n)$-matrix whose entries are holomorphic functions
 (resp.\ holomorphic $1$-forms) on the domain $U$ in $\mathbb{C}^n$,
 \item[(A2)] $B_{\infty}=\operatorname{diag}[\lambda_1,\dots,\lambda_n]$,
where $\lambda_i\in\mathbb{C}$ satisfy $\lambda_i-\lambda_j\not\in \mathbb{Z}\setminus\{0\}$ for $i\not=j$.
\end{enumerate}

Let $H(z)=\det (zI_n-T)$ be the characteristic polynomial of $T$:
\begin{gather*}
 H(z)=\prod_{i=1}^n(z-z_i)=z^n-S_1z^{n-1}+\cdots+(-1)^nS_n,\qquad S_i\in \mathcal{O}_{\mathbb{C}^n}(U).
\end{gather*}
Then we assume the following condition on $H(z)$:
\begin{enumerate}\itemsep=0pt
 \item[(A3)] ${\rm d}S_1\wedge \cdots\wedge {\rm d}S_n\neq 0$ at generic points on~$U$.
\end{enumerate}

It follows from (A3) that the discriminant $\delta_H=\prod_{i<j}(z_i-z_j)^2$ of $H(z)$
is not identically zero and we define a divisor $\Delta_H=\{\delta_H=0\}\cup \{{\rm d}S_1\wedge \cdots\wedge {\rm d}S_n=0\}$ on~$U$.
Taking a simply-connected domain $W\subset U\setminus \Delta_H$ appropriately,
 we can fix a branch of $(z_1,\dots,z_n)$
 and take an invertible matrix~$P$ whose entries are single-valued holomorphic functions on $W$ such that
\begin{equation*}
 P^{-1}TP=\operatorname{diag}[ z_1, \dots, z_n].
\end{equation*}
Decompose $-(zI_n-T)^{-1}B_{\infty}$ into partial fractions
\begin{equation} \label{eq:keibunsum}
 -(zI_n-T)^{-1}B_{\infty}=\sum_{i=1}^n\frac{B_i}{z-z_i}.
\end{equation}
We also assume the condition:
\begin{enumerate}\itemsep=0pt
 \item[(A4)] $r_i:=\operatorname{tr}B_i\neq \pm 1$ on $U\setminus\Delta_H$, $i=1,\dots,n$.
\end{enumerate}

\begin{Remark} At each fixed point $p\in U\setminus\Delta_H$,
 the differential equation in the $z$-direction of~(\ref{eq:okubopfaff}) is an Okubo system,
 which is a Fuchsian system with $n+1$ regular singular points $z=z_1,\dots,z_n,\infty$.
 In the case of $r_i\not\in \mathbb{Z}\setminus \{0\}$, its local monodromy matrix $M_i$ at $z=z_i$ is conjugate to
 \begin{equation*}
 M_i\sim \operatorname{diag}\big[1,\dots,1,e^{2\pi\sqrt{-1}r_i}\big].
 \end{equation*}
 In particular, if $r_i\in\mathbb{Q}\setminus\mathbb{Z}$, the monodromy group is generated by quasi-reflections.
 We use this property in order to introduce flat basic invariants for well-generated finite unitary reflection groups
 in Section~\ref{sec:complexreflection}.
\end{Remark}

\begin{Lemma} \label{lem:lemma2.2}
 Assume that $B_{\infty}$ is invertible.
 Then the Pfaffian system \eqref{eq:okubopfaff} is completely integrable
 if and only if the matrices $T$, $\tilde{\Omega}$, $B_{\infty}$ satisfy the equations
 \begin{gather}
 \big[T,{\tilde \Omega}\big]=O,\qquad {\tilde \Omega}\wedge {\tilde \Omega}=O, \label{eq:commute} \\
 {\rm d}T+\tilde{\Omega}+\big[\tilde{\Omega},B_{\infty}\big]=O, \qquad {\rm d}\tilde{\Omega}=O. \label{eq:ci1}
\end{gather}
 By \eqref{eq:commute}, $T$ and $\tilde{\Omega}$ are simultaneously diagonalizable.
 Indeed, the following holds:
\begin{equation} \label{eq:taikakuBi}
 P^{-1}TP=\operatorname{diag}[ z_1, \dots, z_n ],\qquad
 P^{-1}\tilde{\Omega}P=
 -\operatorname{diag} [ {\rm d}z_1, \dots, {\rm d}z_n ].
\end{equation}
\end{Lemma}

\begin{proof} From the integrability condition ${\rm d}{\rm d}Y=0$ for~(\ref{eq:okubopfaff}),
 we obtain the equations (\ref{eq:ci1}).
 In the following, we shall prove (\ref{eq:taikakuBi}).
 Then the equations (\ref{eq:commute}) follows from (\ref{eq:taikakuBi}).
 Similarly to~(\ref{eq:keibunsum}), we rewrite $-(zI_n-T)^{-1}\tilde{\Omega}B_{\infty}$ into partial fractions
 \begin{equation} \label{eq:keibunsui}
 -(zI_n-T)^{-1}\tilde{\Omega}B_{\infty}=\sum_{i=1}^n\frac{\Omega_i}{z-z_i},
 \end{equation}
 where $\Omega_i$, $1\leq i\leq n$ are $(n\times n)$-matrix valued $1$-forms.
 In order to prove (\ref{eq:taikakuBi}), it is sufficient to show
 \begin{equation} \label{eq:keisurel}
 \Omega_i=-B_i\,{\rm d}z_i,\qquad i=1,\dots,n.
 \end{equation}
 We have
 \begin{equation}
 -\Omega_i-B_i\,{\rm d}z_i+[\Omega_i,B_i]=O,\qquad i=1,\dots,n, \label{eq:ic3okubopfaff}
 \end{equation}
 by substituting (\ref{eq:keibunsum}) and (\ref{eq:keibunsui}) for the integrability condition of~(\ref{eq:okubopfaff}).
 Now let us fix $i\in\{1,\dots,n\}$ arbitrarily.
 We separately consider the cases where $B_i$ is semisimple or not:

1. We assume that $B_i$ is semisimple.
 Note $\operatorname{rank}B_i\leq 1$ and $\operatorname{tr}B_i=r_i\neq \pm 1$ by (A4).
 Then there is a set of linearly independent vectors $\{v_1,\dots,v_n\}\subset \mathbb{C}^n$
 such that $v_j\in \operatorname{ker}B_i$ $(j=1,\dots,n-1)$ and $v_n$ is an eigenvector of $B_i$ belonging to the eigenvalue $r_i$,
 from which we have
 \begin{equation} \label{eq:vecteqm}
 (I_n-B_i)B_i=(1-r_i)B_i.
 \end{equation}
 Noting that $r_i\neq 1$, (\ref{eq:vecteqm}) is equivalent to
 \begin{equation} \label{eq:vecteqm2}
 \bigl((r_i-1)I_n-B_i\bigr)^{-1}B_i=-B_i.
 \end{equation}
 On the other hand, from (\ref{eq:ic3okubopfaff}), we have
 \begin{equation} \label{eq:preeigenv}
(I_n+B_i)\Omega_iv_j=0, \qquad j=1,\dots, n-1, \qquad
 \bigl((r_i-1)I_n-B_i\bigr)\Omega_iv_n={\rm d}z_i\,B_iv_n.
 \end{equation}
 Using (\ref{eq:vecteqm2}), it follows from (\ref{eq:preeigenv}) that
 \begin{equation} \label{eq:eigenOmegai}
 \Omega_iv_j=0, \qquad j=1,\dots,n-1, \qquad
 \Omega_i v_n=-r_i\,{\rm d}z_i\,v_n.
 \end{equation}
 (\ref{eq:eigenOmegai}) implies (\ref{eq:keisurel}).

2. We assume that $B_i$ is not semisimple.
 In this case, $\operatorname{tr} B_i=0$ and $\operatorname{rank}B_i=1$ hold.
 By the assumption, we have $B_i^2=O$.
 In the same manner as the case 1, we see that $\operatorname{ker}B_i\subset \operatorname{ker} \Omega_i$,
 which implies $\Omega_iB_i=O$. From~(\ref{eq:ic3okubopfaff}), we have
 \[
 B_i\Omega_i=B_i\bigl(\Omega_i+B_i\,{\rm d}z_i-[\Omega_i,B_i]\bigr)=O
 \]
 Then we obtain $[\Omega_i,B_i]=O$. Hence (\ref{eq:keisurel}) follows from~(\ref{eq:ic3okubopfaff}).
\end{proof}

\begin{Definition} An {\it extended Okubo system} is a Pfaffian system (\ref{eq:okubopfaff})
 satisfying the system of equations (\ref{eq:commute}) and (\ref{eq:ci1}).
\end{Definition}

\begin{Remark} \label{rm:henkan}
Assume that $T$, $\tilde{\Omega}$, $B_{\infty}$ satisfy the equations $(\ref{eq:commute}),(\ref{eq:ci1})$.
Then replacing $B_{\infty}$ by $B_{\infty}-\lambda I_n$ for any $\lambda\in \mathbb{C}$,
we see that $T$, ${\tilde \Omega}$, $B_{\infty}-\lambda I_n$ also satisfy~(\ref{eq:commute}),~(\ref{eq:ci1}).
 The corresponding transformation from~(\ref{eq:okubopfaff}) into
 \begin{equation*}
 {\rm d}Y^{(\lambda)}=-(zI_n-T)^{-1}\bigl({\rm d}z+\tilde{\Omega}\bigr)(B_{\infty}-\lambda I_n)Y^{(\lambda)}
 \end{equation*}
 is realized by the Euler transformation
 \begin{equation*}
 Y\mapsto Y^{(\lambda)}(z):=\frac{1}{\Gamma(\lambda)}\int (u-z)^{\lambda-1}Y(u)\,{\rm d}u.
 \end{equation*}
\end{Remark}

In virtue of (A3), we can take $(z_1,\dots,z_n)$ as a coordinate system on $W$.
Then, we readily see that (\ref{eq:okubopfaff}) can be rewritten to a Pfaffian system
 \begin{equation} \label{eq:logrep}
 {\rm d}Y=\sum_{i=1}^nB_i\,{\rm d}\log (z-z_i)\,Y.
 \end{equation}
The integrability condition of (\ref{eq:logrep}) is given by
\begin{equation} \label{eq:logrepSch}
 {\rm d}B_i=\sum_{j\neq i}[B_j,B_i]\,{\rm d}\log (z_i-z_j),\qquad i=1,\dots,n,
 \end{equation}
 which is called the {\it Schlesinger system}.
 It is well known (e.g., \cite{JMU,JM}) that isomonodromic deformations of~(\ref{eq:okubo}) are governed
 by the Schlesinger system~(\ref{eq:logrepSch}).
 Therefore we have

 \begin{Lemma} The system of equations \eqref{eq:commute}, \eqref{eq:ci1} is equivalent to the Schlesinger system~\eqref{eq:logrepSch}.
 In particular, if an Okubo system~\eqref{eq:okubo} is given on $\mathbb{P}^1\times \{p_0\}$
 for $p_0\in U\setminus \Delta_H$,
 then there uniquely exists an extended Okubo system on $\mathbb{P}^1\times W$
 which coincides with the given~\eqref{eq:okubo} on $\mathbb{P}^1\times \{p_0\}$,
 where $W\subset U\setminus \Delta_H$ is a neighborhood~of $p_0$.
 \end{Lemma}

 \begin{proof} The lemma follows from the existence and uniqueness of solutions for initial conditions
 of the Schlesinger system shown in~\cite{JMU}.
 \end{proof}

\begin{Remark} \label{Re:z=xn}
Let us assume that $T=T(x)$ in (\ref{eq:okubopfaff}) takes the form
 \begin{equation} \label{cond:T}
 T(x)=-x_n+T_0(x'),\qquad \text{where} \ T_0(x') \ \text{depends only on} \ x'=(x_1,\dots,x_{n-1}),
 \end{equation}
 by taking an appropriate coordinate system $x=(x_1,\dots,x_n)$ on~$U$.
In this case, we may reduce the variable $z$ from the extended Okubo system~(\ref{eq:okubopfaff}):
in~(\ref{eq:okubopfaff}), $z$ and $x_n$ appear in the combination $z+x_n$.
So replacing $z+x_n$ by $x_n$ and omitting ${\rm d}z$,
(\ref{eq:okubopfaff}) becomes
 \begin{equation} \label{eq:okubopfaff0}
 {\rm d}Y_0(x)=T^{-1}\left(\sum_{i=1}^n\tilde{B}^{(i)}\,{\rm d}x_i\right)B_{\infty}Y_0(x),
 \end{equation}
 where $\tilde{B}^{(i)}$ is defined from $\tilde{\Omega}$
 by $\tilde{\Omega}=\sum_{i=1}^n\tilde{B}^{(i)}{\rm d}x_i$.
Then by~(\ref{cond:T}), it holds that
\begin{equation} \label{eq:coeffredOkubo}
 \tilde{B}^{(n)}=I_n, \qquad \tilde{B}^{(i)}=\tilde{B}^{(i)}(x'),\qquad 1\leq i\leq n-1.
\end{equation}
 We call (\ref{eq:okubopfaff0}) the {\it reduced form} of (\ref{eq:okubopfaff}).

 Conversely, given a Pfaffian system (\ref{eq:okubopfaff0}) whose coefficient matrices
 $T$, $\tilde{\Omega}=\sum_{i=1}^n\tilde{B}^{(i)}\,{\rm d}x_i$, $B_{\infty}$ satisfy the conditions
 (\ref{cond:T}), (\ref{eq:coeffredOkubo}), (\ref{eq:commute}), (\ref{eq:ci1}) and the assumptions (A1)--(A4),
 we can recover the extended Okubo system (\ref{eq:okubopfaff}) as a Pfaffian system
 satisfied by
 \[ Y(z,x):=Y_0(x',x_n+z).\]
\end{Remark}

An extended Okubo system naturally induces a Saito bundle in the following manner.
Let us recall the definition of Saito bundle following \cite{DH,Sab}:

\begin{Definition}
 Let $M$ be a complex manifold and $\pi\colon V\to M$ be a holomorphic vector bundle on $M$.
 A {\it Saito bundle} is a $5$-tuple $\big(V, \nabla^V, \Phi^V, R_0^V, R_{\infty}^V\big)$
 consisting of a flat connection $\nabla^V$, a $1$-form $\Phi^V\in\Omega^1(M,\operatorname{End}(V))$
 and two endomorphisms $R_0^V,R_{\infty}^V\in\operatorname{End}(V)$, such that the following conditions are satisfied:
 \begin{gather*}
 \Phi^V\wedge\Phi^V=0,\qquad \big[R_0^V,\Phi^V\big]=0, \\
 {\rm d}^{\nabla^V}\Phi^V=0,\qquad \nabla^V R_0^V+\Phi^V =\big[\Phi^V,R_{\infty}^V\big], \qquad \nabla^V R_{\infty}^V=0,
 \end{gather*}
 where the $\operatorname{End}(V)$-valued forms $\big[R,\Phi^V\big]$
 (with $R:=R_0^V$ or $R_{\infty}^V$), ${\rm d}^{\nabla^V}\Phi^V$ and $\Phi^V\wedge\Phi^V$ are defined by:
 for any $X,Y\in \Gamma(M,\Theta_M)$,
 where $\Theta_M$ denotes the sheaf of vector fields on $M$,
 \begin{gather*}
 \big(\Phi^V\wedge\Phi^V\big)_{X,Y}:=\Phi_X^V\Phi_Y^V-\Phi_Y^V\Phi_X^V,\qquad \big[R,\Phi^V\big]_X:=\big[R,\Phi_X^V\big], \\
 \big({\rm d}^{\nabla^V}\Phi\big)_{X,Y}:=\nabla_X^V\big(\Phi_Y^V\big)-\nabla^V_{Y}\big(\Phi_X^V\big)-\Phi_{[X,Y]}^V.
 \end{gather*}
 \end{Definition}

Returning to our setting, let $U\times \mathbb{C}^n$ be a trivial bundle of rank $n$ on the domain $U\subset \mathbb{C}^n$
and $({\bf e}_1,\dots,{\bf e}_n)$ be a basis of the fiber space $\mathbb{C}^n$.
We identify the matrices $\tilde{\Omega}^t$, $T^t$, $B_{\infty}^t=B_{\infty}$ as endomorphisms of $U\times \mathbb{C}^n$ respectively by
\[
 \tilde{\Omega}^t({\bf e}_j):=\sum_{i=1}^n\tilde{\Omega}_{ji}{\bf e}_i,\qquad
 T^t({\bf e}_j):=\sum_{i=1}^nT_{ji}{\bf e}_i,\qquad
 B_{\infty}({\bf e}_j):=\sum_{i=1}^n(B_{\infty})_{ji}{\bf e}_i
\]
Here, for a square matrix $A$, $A^t$ denotes the transpose of $A$.
Then, for an extended Okubo system~(\ref{eq:okubopfaff}),
Lemma~\ref{lem:lemma2.2} implies that the $5$-tuple $\big(U\times \mathbb{C}^n,{\rm d},\alpha\tilde{\Omega}^t,\alpha T^t,B_{\infty}\big)$ with some $\alpha\in\mathbb{C}\setminus \{0\}$ defines a~Saito bundle on~$U$.

\section{Logarithmic vector fields and free divisor} \label{sec:logvecf}

In this section, we study a divisor defined by the singular locus of an extended Okubo system.
The main purpose of this section is to prove that the set of the vector fields defined by~(\ref{equ:LogVFs})
forms a unique standard system of generators of logarithmic vector fields along the divisor
when the divisor is free (Lemma~\ref{lemms:non-negative}, cf.\ Remark~\ref{rem:flatvector}).

\subsection{Logarithmic vector fields}

Let a domain $U\subset \mathbb{C}^n$ have the form of
$U=U'\times \mathbb{C}\subset \mathbb{C}^{n-1}\times\mathbb{C}$.
We consider a monic polynomial $h(x)=h(x',x_n)$ in $x_n$ of degree $n$:
\begin{equation} \label{poly:h}
 h(x)=x_n^n-s_1(x')x_n^{n-1}+\dots+(-1)^n s_n(x')
 =\prod_{i=1}^n\big(x_n-z_i^0(x')\big),
\end{equation}
where $s_i(x')\in \mathcal{O}_{\mathbb{C}^{n-1}}(U')$.

Let $D$ be the divisor on $U$ defined by $h(x)$,
i.e., $D=\{(x)\in U;\, h(x)=0\}$. We give definitions of logarithmic vector field along $D$ and free divisor following K.~Saito~\cite{Sai1}:

\begin{Definition}
Let $\mathcal{M}_{\mathbb{C}^{n-1}}(U')$ be the field of meromorphic functions on $U'$.
A vector field
$V=\sum_{k=1}^n v_k\partial_{x_k}$
with $v_k\in \mathcal{M}_{\mathbb{C}^{n-1}}(U')\otimes_{\mathbb{C}}\mathbb{C}[x_n]$ is called a meromorphic
logarithmic vector field along $D$
if $(V h)/h\in \mathcal{M}_{\mathbb{C}^{n-1}}(U')\otimes_{\mathbb{C}}\mathbb{C}[x_n]$,
or equivalently if $(V h)|_D=0$.
If moreover, $v_k\in \mathcal{O}_{\mathbb{C}^{n-1}}(U')\otimes_{\mathbb{C}}\mathbb{C}[x_n]$,
then $V$ is called a logarithmic vector field along $D$.

Let $\operatorname{Der}(-\log D)$ be the set of logarithmic vector fields along $D$,
 which forms naturally an $\mathcal{O}_{\mathbb{C}^{n-1}}(U')\otimes _{\mathbb{C}}\mathbb{C}[x_n]$-module.
The divisor $D$ is said to be free if
$\operatorname{Der}(-\log D)$ is a free $\mathcal{O}_{\mathbb{C}^{n-1}}(U')\otimes_{\mathbb{C}} \mathbb{C}[x_n]$-module.
\end{Definition}

For the set of meromorphic logarithmic vector fields
$\mathcal{V}=\big\{V_i=\sum_{j=1}^n v_{ij}\partial_{x_j}\big\}_{1\leq i\leq n}$ along $D$,
let $M_{\mathcal{V}}(x)$ denote the $(n\times n)$-matrix whose
$(i,j)$-entry is defined by $v_{n-i+1,j}(x)$.

\begin{Remark} \label{rem:Saitocriterion}
It is known that $D$ is free if there is
$\mathcal{V}=\{V_i\}_{1\leq i\leq n}\subset\operatorname{Der}(-\log D)$ such that
$\det M_{\mathcal{V}}=h$ (Saito's criterion \cite{Sai1}).
In this case, the matrix $M_{\mathcal{V}}$ is called a {\it Saito matrix}.
\end{Remark}

Now introduce holomorphic functions $S_1(x),\dots,S_n(x)$ on $U$ by
\begin{equation} \label{eq:H(z,x)fromh(x)}
 H(z,x):=h(x',x_n+z)=z^n-S_1(x)z^{n-1}+\cdots +(-1)^nS_n(x).
\end{equation}
We rewrite the condition (A3) for (\ref{eq:H(z,x)fromh(x)}) in terms of a condition for meromorphic logarithmic vector fields:

\begin{Lemma} \label{lemma:zeroV}
The following conditions $(i)$ and $(ii)$ are equivalent.
\begin{enumerate}\itemsep=0pt
\item[$(i)$] ${\rm d}S_1\wedge \cdots \wedge {\rm d}S_n\neq 0$ on generic points of~$U$ $($the condition~{\rm (A3))}.
\item[$(ii)$] Let $V$ be any meromorphic logarithmic vector field along $D$ satisfying $[\partial_{x_n},V]=0$. Then $V=0$.
\end{enumerate}
\end{Lemma}

\begin{proof}
Any meromorphic vector field $V$ satisfying $[\partial_{x_n},V]=0$ is written as $V\!=\!\sum_{k=1}^n \! v_k(x')\partial_{x_k}$, $v_k(x')\in \mathcal{M}_{\mathbb{C}^{n-1}}(U')$.
Then
\begin{gather*}
 Vh= \left(nv_n-\sum_{k=1}^{n-1}v_k\frac{\partial s_1}{\partial x_k}\right)x_n^{n-1}+\left(-v_n(n-1)s_1
 +\sum_{k=1}^{n-1}v_k\frac{\partial s_2}{\partial x_k}\right)x_n^{n-2} \\
\hphantom{Vh=}{} +\cdots +(-1)^n\left(-v_ns_{n-1}+\sum_{k=1}^{n-1}v_k\frac{\partial s_n}{\partial x_k}\right).
 \end{gather*}
 So the condition $Vh=0$ is equivalent to
\begin{equation}\label{equ:zeroV-2}
 \begin{pmatrix}
 v_1(x')&\dots&v_n(x')
 \end{pmatrix}
 \begin{pmatrix}
 \dfrac{\partial s_1(x')}{\partial x_1}&\dfrac{\partial s_2(x')}{\partial x_1}
 &\dots&\dfrac{\partial s_n(x')}{\partial x_1}\\
 \dots\\
 \dfrac{\partial s_1(x')}{\partial x_{n-1}}&
 \dfrac{\partial s_2(x')}{\partial x_{n-1}}
 &\dots&\dfrac{\partial s_n(x')}{\partial x_{n-1}}\\
 -n&-(n-1)s_1(x')&\dots&-s_{n-1}(x')
 \end{pmatrix}
 =O.
\end{equation}
The linear equation (\ref{equ:zeroV-2}) has a non-zero solution vector $(v_1(x'),\dots,v_n(x'))$ if and only if
${\rm d}S_1\wedge \cdots \wedge {\rm d}S_n\equiv 0$ holds for $H(z,x)$ in~(\ref{eq:H(z,x)fromh(x)}). This proves the lemma.
\end{proof}

Let us take a single-valued branch of $z_i^0(x')$ $(i=1,\dots,n)$ in~(\ref{poly:h}) on an appropriate domain $W'\subset U'$.
Then we define matrices $P_h(x'), M_h(x'), M_{\mathcal{V}^{(h)}}(x)$ by
 \begin{gather}\label{equ:Def of P_h}
 P_h(x')=
 \left(\frac{\partial z_j}{\partial x_i}\right),\\
\label{equ:Def of M_V and M_h}
 M_{h}(x')=P_h \operatorname{diag}\big[z_1^0,z_2^0,\dots,z_n^0\big] P_h^{-1},\qquad
 M_{\mathcal{V}^{(h)}}(x)=-P_h \operatorname{diag}[z_1,z_2,\dots,z_n] P_h^{-1},
\end{gather}
and vector fields $V_i^{(h)}$, $1\leq i \leq n$ by
\begin{equation}\label{equ:LogVFs}
 \big(\begin{matrix}
 V^{(h)}_n & \cdots & V^{(h)}_1
 \end{matrix}\big)^t
 = M_{\mathcal{V}^{(h)}}
 \begin{pmatrix}
 \partial_{x_1} & \cdots & \partial_{x_n}
 \end{pmatrix}^t.
\end{equation}
By definition, we have
\begin{equation}\label{equ:rel of M_V and M_h}
 M_{\mathcal{V}^{(h)}}(x)= x_n I_n-M_{h}(x')
\end{equation}
and
\begin{equation}\label{equ:det(x_nI_n-M_h(x'))}
 \det M_{\mathcal{V}^{(h)}}(x)=h(x).
\end{equation}

\begin{Lemma}\label{lemma:LogVFs}
The vector fields $V^{(h)}_i$, $1\le i\le n$ are meromorphic logarithmic vector fields along~$D$.
\end{Lemma}

\begin{proof} First we prove that the entries of $M_h(x')$ are meromorphic functions on $U'$. Put
 \[
 F=
 \begin{pmatrix} (z_1^0)^{n-1}&\dots&(z_n^0)^{n-1} \\
 \vdots & \ddots & \vdots \\
 z_1^0&\dots&z_n^0\\1&\dots&1
 \end{pmatrix}.
 \]
 Introduce two matrices $M^0$ and $M^1$ respectively by
 \begin{equation*}
 M^0=F\,\operatorname{diag}\big[z_1^0,z_2^0,\dots,z_n^0\big]F^{-1}\qquad \text{and}\qquad
 M^1=P_hF^{-1}.
 \end{equation*}
 Then the $(i,j)$-entry of $M^0$ is given by
 \[
 (M^0)_{ij}=\frac{1}{\det F}
 \begin{vmatrix}
 \big(z_1^0\big)^{n-1} & \cdots & \big(z_n^0\big)^{n-1} \\
 \vdots & \cdots & \vdots \\
 \big(z_1^0\big)^{n+1-i} & \cdots & \big(z_n^0\big)^{n+1-i} \\
 \vdots & \cdots & \vdots \\
 1 & \cdots & 1
 \end{vmatrix}\text{$\leftarrow$ $j$-th row},
 \]
 which is symmetric with respect to permutations of~$z_i^0$. Therefore $M^0$ is meromorphic on~$U'$. We can see in a similar way that~$M^1$ is also meromorphic on $U'$. Since $M_h$ is written as
 \begin{equation*}
 M_h=M^1M^0\big(M^1\big)^{-1},
 \end{equation*}
we see that the entries of $M_h(x')$ are meromorphic functions on~$U'$.
 Then it is clear that $V_i^{(h)}$, $1\leq i\leq n$ are meromorphic logarithmic vector fields
 from $z_i\partial_{z_i}h=h$ $(i=1,\dots,n)$ and~(\ref{equ:Def of P_h}), (\ref{equ:Def of M_V and M_h}), (\ref{equ:LogVFs}).
\end{proof}

\begin{Lemma} \label{lem:uniqunesslogvf}
Assume that there are $V_i\in\operatorname{Der}(-\log D)$, $1\leq i\leq n$, satisfying $[\partial_{x_n},V_i-x_n\partial_{x_n+1-i}]=0$.
Then it holds that $V_i=V_i^{(h)}$, $1\le i\le n$, which implies $\mathcal{V}^{(h)}=\big\{V_i^{(h)}\big\}_{1\leq i\leq n}\subset \operatorname{Der}(-\log D)$, and hence $D$ is free.
\end{Lemma}

\begin{proof} By (\ref{equ:rel of M_V and M_h}) and (\ref{equ:LogVFs}), we see $\big[\partial_{x_n},V_i-V_i^{(h)}\big]=0$,
 which implies $V_i-V_i^{(h)}=0$ in virtue of Lemma~\ref{lemma:zeroV}. Then the equality~(\ref{equ:det(x_nI_n-M_h(x'))}) implies that $M_{\mathcal{V}^{(h)}}$ is a Saito matrix and $D$ is free.
\end{proof}

\subsection{Quasi-homogeneous polynomial case}\label{sec:weighted h}

In this subsection we assume that the function $h(x)=h(x',x_n)$ in (\ref{poly:h}) is a polynomial in $x$ and weighted
homogeneous with respect to a weight
$w(\cdot)$ with
\[
 0<w(x_1)\le w(x_2)\le\dots \le w(x_{n-1})\le w(x_n)=1.
\]
In this case, we replace $\mathcal{O}_{\mathbb{C}^{n-1}}(U')$ and $\mathcal{M}_{\mathbb{C}^{n-1}}(U')$
in the previous subsection by $\mathbb{C}[x']$ and $\mathbb{C}(x')$ respectively.
\begin{Lemma}\label{lemma:V1=Eu}
 It holds that $V^{(h)}_1=\sum_{k=1}^n w(x_k)x_k \partial_{x_k}$, namely, $V^{(h)}_1$
 is an Euler operator.
\end{Lemma}

\begin{proof}
 The lemma is clear from that the Euler operator is a logarithmic vector field along $D$
and the argument similar to the proof of Lemma~\ref{lem:uniqunesslogvf}.
\end{proof}

Note $w\big(z_m^0(x')\big)=w(x_n)=1$.
From the definition of $M_h$ in (\ref{equ:Def of M_V and M_h}), we find
$w((M_h)_{i,j})
=w\big(z_m^0\big)+w\big({\partial z_m^0\over \partial x_i}\big)-w\big({\partial z_m^0\over \partial x_j}\big)$ for $1\le i,j\le n$, that is,
\begin{equation*}
 w ((M_h)_{i,j} )=1-w(x_i)+w(x_j),\qquad 1\le i,j\le n.
\end{equation*}
Since $w ((M_h)_{i,j} )=w\big(\big(M_{\mathcal{V}^{(h)}}\big)_{i,j}\big)$,
we have a kind of duality
\begin{equation*}
 w(x_i)+w\big(V^{(h)}_{n-i+1}\big)=1,\qquad 1\le i\le n,
\end{equation*}
between $\{w(x_i)\}$ and $\big\{w\big(V_i^{(h)}\big)\big\}$.

\begin{Lemma}\label{lemms:wVige0}
Let $V\in\operatorname{Der}(-\log D)$ be weighted homogeneous.
 Then $w(V)\ge 0$.
\end{Lemma}

\begin{proof}
Let $V=\sum_{j=1}^n v_j\partial_{x_j}$.
If $w(V)<0$, then $v_{j}\in \mathbb{C}[x']$.
This implies $V=0$ by Lem\-ma~\ref{lemma:zeroV}.
\end{proof}

\begin{Lemma}\label{lemms:generatorsVk}
Assume that there exists a set $\mathcal{V}=\{V_i\}_{1\leq i\leq n}\subset\operatorname{Der}(-\log D)$
satisfying
\[ \det M_{\mathcal{V}}(x)|_{x'=0}=cx_n^n\]
for $c\in\mathbb{C}\setminus \{0\}$
and $V_i$, $1\leq i\leq n$ are weighted homogeneous with $w(V_1)\le w(V_2)\le\dots\le w(V_n)$.
Then the following three assertions hold.

\begin{enumerate}\itemsep=0pt
\item[$(i)$]
$\det M_{\mathcal{V}}(x)=ch(x)$ and thus $D$ is free.

\item[$(ii)$]
$w(V_i)=w\big(V^{(h)}_i\big)$ $(=1-w(x_{n-i+1}))$, $1\le i\le n$.

\item[$(iii)$]
 There is a matrix $G(x')\in {\rm GL}(n,\mathbb{C}[x'])$ such that
 $w(G(x')_{i,j})=w(x_j)-w(x_i)$ and
 \[
 \big(\begin{matrix}
 V^{(h)}_n&V^{(h)}_{n-1}& \dots& V^{(h)}_1
 \end{matrix}\big)^t
 =G(x')
 \begin{pmatrix}
 V_n&V_{n-1}& \dots& V_1
 \end{pmatrix}^t,
 \]
which implies $\mathcal{V}^{(h)}=\big\{V_i^{(h)}(x)\big\}_{1\le i\le n}\subset\operatorname{Der}(-\log D)$.
\end{enumerate}
\end{Lemma}

\begin{proof}
 Proof of~(i). Since $V_i$, $1\leq i\leq n$ are logarithmic vector fields along $D$,
 it holds that $\det M_{\mathcal{V}}=f(x)h(x)$ for some $f(x)\in \mathbb{C}[x]$ \cite[(1.5)iii)]{Sai1}.
 Then it follows that $\det M_{\mathcal{V}}=ch$.

 Proof of (ii).
 Since $w(V_i)\ge 0$ from Lemma~\ref{lemms:wVige0},
 it holds that $(M_{\mathcal{V}})|_{x=0}=0$.
 Then since
 $\det(M_{\mathcal{V}})|_{x'=0}=c x_n^n$,
it holds that
\begin{equation*}
 M_{\mathcal{V}}(x)|_{x'=0}=x_n R,
\end{equation*}
where $R\in {\rm GL}(n,\mathbb{C})$.
This proves $w(V_i)=1-w(x_{n-i+1})$.

Proof of (iii).
Note that $R^{-1} M_{\mathcal{V}}(x)|_{x'=0}=x_n I_n.$
By the argument above, we see that $R^{-1} M_{\mathcal{V}}(x)$ is written as
\begin{equation*}
 R^{-1} M_{\mathcal{V}}(x)= x_n \tilde{R}(x')+\tilde{\tilde{R}}(x'),
\end{equation*}
for some $\tilde{R}(x'),\tilde{\tilde{R}}(x')\in \mathbb{C}[x']^{n\times n}$
satisfying
 \[
 \tilde{R}(0)=I_n,\qquad w(\tilde{R}_{i,j})=-w(x_i)+w(x_j),\qquad
 w(\tilde{\tilde{R}}_{i,j})=1-w(x_i)+w(x_j).
 \]
Then $\det \tilde{R}(x')=1$ and
$\tilde{R}(x')^{-1}R^{-1}M_{\mathcal{V}}-x_nI_n\in \mathbb{C}[x']^{n\times n}$,
which (together with Lemma \ref{lem:uniqunesslogvf}) proves
 \[
 \tilde{R}(x')^{-1}R^{-1}
 \begin{pmatrix}
 V_n&V_{n-1}& \dots& V_1
 \end{pmatrix}^t
 =
 \big(\begin{matrix}
 V^{(h)}_n&V^{(h)}_{n-1}& \dots& V^{(h)}_1
 \end{matrix}\big)^t.\tag*{\qed}
 \]\renewcommand{\qed}{}
\end{proof}

\begin{Lemma}\label{lemms:non-negative}
 Assume that $D$ is free.
 Then the following assertions hold:
 \begin{enumerate}\itemsep=0pt
 \item[$(i)$]
 $\mathcal{V}^{(h)}=\{V^{(h)}_i\}_{1\le i\le n}$ is a system of generators of $\operatorname{Der}(-\log D)$
 satisfying $M_{\mathcal{V}^{(h)}}-x_nI_n\in \mathbb{C}[x']^{n\times n}$.
 \item[$(ii)$]
 Let $\mathcal{V}=\{V_i\}_{1\le i\le n}$ be any system of generators of $\operatorname{Der}(-\log D)$
 satisfying $M_{\mathcal{V}}-x_nI_n\in\mathbb{C}[x']^{n\times n}$.
 Then $\mathcal{V}=\mathcal{V}^{(h)}$.
 \end{enumerate}
\end{Lemma}

\begin{proof}
We prove (i).
By the assumption that $D$ is free,
there is a set of generators $\mathcal{V}=\{V_i\}_{1\leq i\leq n}\subset\operatorname{Der}(-\log D)$ such that $\det(M_{\mathcal{V}})=h(x)$.
Since each homogeneous part of $V_i$ are logarithmic vector field
with a non-negative weight (Lemma~\ref{lemms:wVige0}),
all entries of $(M_{\mathcal{V}})|_{x'=0}$ are polynomials in $x_n$ of positive
degrees.
Put $R(x_n):=x_n^{-1}(M_{\mathcal{V}})|_{x'=0}\in \mathbb{C}[x_n]^{n\times n}$.
Then $\det R(x_n)=1$, and hence $\det R(0)=1$.
Let $V'_i$, $1\le i\le n$ be logarithmic vector fields defined by
\[
\begin{pmatrix} V'_n&V'_{n-1}& \dots& V'_1 \end{pmatrix}^t
=R(0)^{-1}M_{\mathcal{V}}
\begin{pmatrix} \partial_{x_1}&\partial_{x_2}& \dots&
\partial_{x_n} \end{pmatrix}^t.
\]
Put
$M_{\mathcal{V}'}=R(0)^{-1}M_{\mathcal{V}}$.
Then $(M_{\mathcal{V}'})|_{x'=0}=x_n I_n+O\big(x_n^2\big)$.
Let $V''_i$ be the homogeneous part of~$V'_i$ with
$w(V''_i)=1-w(x_{n-i+1})$.
Then it holds that $\det (M_{\mathcal{V}''})|_{x'=0}=x_n^n$, and hence
Lemma~\ref{lemms:generatorsVk} implies that
$V^{(h)}_i$, $1\le i\le n$ are logarithmic vector fields.
(ii) follows from Lemma~\ref{lem:uniqunesslogvf}.
\end{proof}

\begin{Remark} \label{rem:flatvector}
 In the case where $h(x)$ is the discriminant of a well-generated complex reflection group,
 D.~Bessis \cite{Be} showed the existence of a system of generators of $\operatorname{Der}(-\log D)$
 whose Saito matrix $M_{\mathcal{V}}$ has the form $M_{\mathcal{V}}-x_nI_n\in\mathbb{C}[x']^{n\times n}$.
 Such a system of generators is called {\it flat} in~\cite{Be}.
 This can be verified also by explicit representations of basic derivations of well-generated complex reflection groups
 found in~\cite{BM,Or, OT}.
\end{Remark}

\section{Saito structure (without metric)} \label{sec:saitostrandOkubo}

\subsection{Review on Saito structure (without metric)} \label{subsec:Saito}

In this subsection, we review Saito structure (without metric) introduced by Sabbah~\cite{Sab}.
(In the sequel, we abbreviate Saito structure (without metric) to Saito structure for brevity.)

\begin{Definition}[C.~Sabbah~\cite{Sab}] \label{def:Saito}
 Let $M$ be a complex manifold of dimension $n$,
 $TM$ be its tangent bundle,
 and $\Theta_M$ be the sheaf of holomorphic sections of $TM$.
 A {\it Saito structure} on $M$ is a data consisting of $(\knabla, \Phi, e, E)$ in (i)-(iii)
 satisfying the conditions (a), (b) below:
 \begin{enumerate}\itemsep=0pt
 \item[(i)] $\knabla$ is a flat torsion-free connection on $TM$,
 \item[(ii)] $\Phi$ is a symmetric Higgs field on $TM$,
 \item[(iii)] $e$ and $E$ are vector fields (i.e., global sections of $TM$),
 respectively called {\it unit field} and {\it Euler field}.
 \end{enumerate}

 \begin{enumerate}\itemsep=0pt
 \item[(a)] A meromorphic connection $\nabla$ on the bundle $\pi^{*}TM$ on $\mathbb{P}^1\times M$ defined by
 \begin{equation} \label{def:nabla}
 \nabla=\pi^{*}\knabla+\frac{\pi^*\Phi}{z}-\left(\frac{\Phi(E)}{z}+\knabla E\right)\frac{dz}{z}
 \end{equation}
 is integrable, where $\pi$ is the projection $\pi :\mathbb{P}^1\times M\rightarrow M$
 and $z$ is a non-homogeneous coordinate of $\mathbb{P}^1$,
 \item[(b)] the field $e$ is $\knabla$-horizontal (i.e., $\knabla (e)=0$) and satisfies $\Phi_e={\rm Id}$,
 where we regard $\Phi$ as an $\operatorname{End}_{\mathcal{O}_M}(\Theta_M)$-valued $1$-form
 and $\Phi_e\in \operatorname{End}_{\mathcal{O}_M}(\Theta_M)$ denotes the contraction of the vector field $e$ and the $1$-form $\Phi$.
 \end{enumerate}
\end{Definition}

\begin{Remark} To the Higgs field $\Phi$ there associates a product $\star$ on $\Theta_M$ defined by
$X \star Y =\Phi_{X}(Y)$ for $X,Y\in \Theta_M$.
 The Higgs field $\Phi$ is said to be symmetric if the product $\star$ is commutative and associative.
 The condition $\Phi_e={\rm Id}$ in Definition \ref{def:Saito}(b) implies that the field $e$ is the unit of the product $\star$.
 So we can introduce on $\Theta_M$ the structure of associative and commutative $\mathcal{O}_M$-algebra with a unit.
\end{Remark}

Since the connection $\knabla$ is flat and torsion free, we can take a {\it flat coordinate system} $(t_1,\dots,t_n)$
such that $\knabla(\partial _{t_i})=0$ $(i=1,\dots,n)$
at least on a~simply-connected open set~$U$ of~$M$.
We take a~flat coordinate system $(t_1,\dots,t_n)$ on~$U$
and assume the following conditions:
\begin{enumerate}\itemsep=0pt
 \item[(C1)] $e=\partial_{t_n}$,
 \item[(C2)] $E=w_1t_1\partial _{t_1}+\cdots+w_nt_n\partial _{t_n}$ for $w_i\in\mathbb{C}$ $(i=1,\dots,n)$,
 \item[(C3)] $w_n=1$ and $w_i-w_j\not\in\mathbb{Z}\setminus\{0\}$ for $i\neq j$.
\end{enumerate}
In this paper, a function $f\in\mathcal{O}_M(U)$ is said to be weighted homogeneous with a weight $w(f)\in \mathbb{C}$
if $f$ is an eigenfunction of the Euler operator: $Ef=w(f)f$.
In particular, the flat coordina\-tes~$t_i$, $1\leq i \leq n$ are weighted homogeneous with $w(t_i)=w_i$.

We fix a basis $\{\partial _{t_1},\dots,\partial _{t_n}\}$ of $\Theta_M(U)$
and introduce the following matrices:
\begin{enumerate}\itemsep=0pt
 \item[(i)] $\tilde{\Phi}=(\tilde{\Phi}_{ij})$ is (the transpose of) the representation matrix of $\Phi$,
 namely the $(i,j)$-entry $\tilde{\Phi}_{ij}$ is a $1$-form on $U$ defined by
 \begin{equation*}
 \Phi(\partial _{t_j})=\sum_{i=1}^n\tilde{\Phi}_{ji}\partial _{t_i},
\qquad i=1,\ldots,n,
 \end{equation*}
 \item[(ii)] $\mathcal{T}=(\mathcal{T}_{ij})$ and $\mathcal{B}_{\infty}=\bigl((\mathcal{B}_{\infty})_{ij}\bigr)$
 are (the transpose of) the representation matrices of $-\Phi(E)$ and $\knabla E$ respectively,
 namely
 \begin{equation*}
 -\Phi_{\partial_{t_j}}(E)=\sum_{i=1}^n\mathcal{T}_{ji}\partial_{t_i},\qquad
 \knabla_{\partial_{t_j}}(E)=\sum_{i=1}^n(\mathcal{B}_{\infty})_{ji}\partial_{t_i}.
 \end{equation*}
\end{enumerate}

We assume that $-\Phi(E)$ is generically regular semisimple on $M$,
that is the discriminant of $\det(z-\mathcal{T})$ does not identically vanish on $M$.

\begin{Lemma} \label{lem:Btow}
The matrix $\mathcal{B}_{\infty}$ is a constant diagonal matrix: $\mathcal{B}_{\infty}=\operatorname{diag}[w_1,\dots,w_n]$.
\end{Lemma}

\begin{proof} It is straightforward.
\end{proof}

\begin{Lemma} \label{lem:trivial}
 The meromorphic connection $\nabla$ in $(\mbox{{\rm a}})$ is integrable if and only if $\mathcal{T}$, $\tilde{\Phi}$
 and $\mathcal{B}_{\infty}$ are subject to the following relations
 \begin{gather}
 \bigl[\mathcal{T},\tilde{\Phi}\bigr]=O, \quad
 \tilde{\Phi}\wedge\tilde{\Phi}=O, \label{eq:saitorel1} \\
 {\rm d}\mathcal{T}+\tilde{\Phi}+\big[\tilde{\Phi},\mathcal{B}_{\infty}\big]=O, \qquad d\tilde{\Phi}=O. \label{eq:saitorel2}
 \end{gather}
\end{Lemma}

\begin{proof}
 See \cite[equation~(2.6), p.~203]{Sab}.
\end{proof}

Lemma~\ref{lem:trivial} implies that a Saito structure induces an extended Okubo system
and thus a~Saito bundle $\big(U\times \mathbb{C}^n,{\rm d},\alpha\tilde{\Phi}^t,\alpha\mathcal{T}^t,\mathcal{B}_{\infty}\big)$
with some $\alpha\in\mathbb{C}\setminus \{0\}$
by taking a flat coordinate system on~$U$.
In the next subsection, we show the opposite direction,
namely the space of independent variables of an extended Okubo system satisfying some generic condition
can be equipped with a Saito structure.

 \begin{Remark} \label{rem:BirkhoffOkubo}
 The meromorphic connection (\ref{def:nabla}) is written in the following matrix form with respect to the flat coordinate system:
 \begin{equation} \label{eq:Birkhoff}
 {\rm d} \mathcal{Y}=\left(\left(-\frac{\mathcal{T}}{z}+\mathcal{B}_{\infty}\right) \frac{{\rm d}z}{z} -\frac{1}{z}\tilde{\Phi}\right)\mathcal{Y}.
 \end{equation}
 The system of equations (\ref{eq:saitorel1}),(\ref{eq:saitorel2}) is equivalent to the integrability condition of (\ref{eq:Birkhoff}).
 The system of ordinary linear differential equations
 \begin{equation*}
 \frac{{\rm d}\mathcal{Y}}{{\rm d}z}=\left(-\frac{\mathcal{T}}{z^2}+\frac{\mathcal{B}_{\infty}}{z}\right)\mathcal{Y}
 \end{equation*}
 has an irregular singularity of Poincar\'e rank one at $z=0$ and a regular singularity at $z=\infty$,
 which is called a {\it Birkhoff normal form}.
 A Birkhoff normal form can be transformed into an Okubo system using a Fourier--Laplace transformation.
 \end{Remark}

\begin{Lemma} \label{lem:Tlogvf}
 Define vector fields $V_i$, $1\leq i\leq n$ by
 \begin{equation} \label{eq:Saito logvec}
 \begin{pmatrix} V_n & \cdots & V_1 \end{pmatrix}^t=-\mathcal{T}\begin{pmatrix} \partial_{t_1} & \cdots & \partial_{t_n} \end{pmatrix}^t
 \end{equation}
 and put $h=h(t)=\det (-\mathcal{T})$.
 Then
 $V_i$, $1\leq i\leq n$ are logarithmic vector fields along $D=\{t\in X;$ $h(t)=0\}$,
 and $D$ is a free divisor.
\end{Lemma}

\begin{proof} Several proofs are found in \cite{KMS, Sab}.
\end{proof}

\begin{Lemma} \label{lem:Cishomo}
 There is a unique matrix $\mathcal{C}$ whose entries are in $\mathcal{O}_M(U)$ such that
 \begin{equation*}
 \mathcal{T}=-E\mathcal{C}, \qquad
 \tilde{\Phi}={\rm d}\mathcal{C}
 \end{equation*}
 and that each matrix entry $\mathcal{C}_{ij}$ of $\mathcal{C}$ is weighted homogeneous with $w(\mathcal{C}_{ij})=1-w_i+w_j$.
\end{Lemma}

\begin{proof} The lemma follows from (\ref{eq:saitorel2}).
\end{proof}

\begin{Lemma} \label{lem:flatC=t}
 Let $(t_1,\dots,t_n)$ be a flat coordinate system of a Saito structure.
 Then it holds that $t_j=\mathcal{C}_{nj}=-w_j^{-1}\mathcal{T}_{nj}$, $1\leq j\leq n$.
\end{Lemma}

\begin{proof}
Let $V_i$, $1\leq i\leq n$ be same as (\ref{eq:Saito logvec}).
Lemma~\ref{lem:Tlogvf} shows that each $V_i$ is
logarithmic along~$D$.
 Then similarly to Lemma~\ref{lemma:V1=Eu}, it holds that $V_1=E$,
 which implies $-\mathcal{T}_{nj}=w_jt_j$.
\end{proof}

\begin{Proposition}[Manin \cite{Ma}, Konishi--Minabe \cite{KM}] \label{prop:potential vector field} There is a unique $n$-tuple of holomorphic functions $\vec{g}=(g_1,\dots,g_n)\in \mathcal{O}_M^n(U) $ such that $\mathcal{C}_{ij}=\frac{\partial g_j}{\partial t_i}$, and that $g_j$ is weighted homogeneous with $w(g_j)=1+w_j$. The vector $\vec{g}$ $($or precisely the vector field $\mathcal{G}=\sum_{i=1}^ng_i\partial_{t_i})$ is called a~potential vector field.
\end{Proposition}

\begin{Remark} \label{rem:pvfC} It is readily found that the potential vector field $\vec{g}=(g_1,\dots,g_n)$ is explicitly given by
 \begin{equation} \label{eq:pvfC}
 g_j=\frac{1}{1+w_j}\sum_{i=1}^nw_it_i\mathcal{C}_{ij},\qquad 1\leq j\leq n.
 \end{equation}
\end{Remark}

\begin{Remark} \label{rem:flat pvf constant}
 By definition, the flat coordinate $t=(t_1,\dots,t_n)$ has ambiguity of weight preserving linear transformations.
 In particular, if $t_i$, $1\le i\le n$ is changed to $t_i'=c_it_i$ for non-zero constants $c_i$,
 the corresponding potential vector field is changed to $g_i'(t')=c_ic_ng_i(t)$.
\end{Remark}

\begin{Proposition} \label{prop:genWDVV}
 The potential vector field $\vec{g}=(g_1,\dots,g_n)$ is a solution to the following system of nonlinear equations:
 \begin{gather}
 \sum_{m=1}^n\frac{\partial^2 g_m}{\partial t_k\partial t_i}\frac{\partial^2 g_j}{\partial t_l\partial t_m}=
 \sum_{m=1}^n\frac{\partial^2 g_m}{\partial t_l\partial t_i}\frac{\partial^2 g_j}{\partial t_k\partial t_m},
\qquad i,j,k,l=1,\dots,n, \label{eq:genWDVV} \\
 \frac{\partial ^2 g_j}{\partial t_n\partial t_i}=\delta_{ij},\qquad i,j=1,\dots,n, \label{eq:vpunit} \\
 Eg_j=\sum_{k=1}^nw_kt_k\frac{\partial g_j}{\partial t_k}=(1+w_j)g_j,\qquad j=1,\dots,n. \label{eq:vphomo}
 \end{gather}
\end{Proposition}

\begin{proof}
 The equation~(\ref{eq:genWDVV}) follows from the associativity of $\star$ (i.e., $\tilde{\Phi}\wedge\tilde{\Phi}=0$).
 The equation~(\ref{eq:vpunit}) follows from $\frac{\partial\mathcal{C}}{\partial t_n}=I_n$ (i.e., $\Phi_e=\mbox{Id}$).
\end{proof}

\begin{Definition} \label{def:extWDVV}
The system of non-linear differential equations
(\ref{eq:genWDVV})--(\ref{eq:vphomo})
for a vector $\vec{g}=(g_1,\ldots,g_n)$ is called the extended WDVV equation.
\end{Definition}

\begin{Remark}
 The notion of ``$F$-manifolds with compatible flat structures" was introduced by Manin \cite{Ma} as a generalization of Frobenius manifolds. ``Potential vector field'' in Proposition~\ref{prop:potential vector field} is called ``local vector potential'' in Manin's framework~\cite{Ma}.
 And the associativity conditions~(\ref{eq:genWDVV}),~(\ref{eq:vpunit}) are called ``oriented associativity equations"
 in~\cite{Ma}.
 The authors were informed these facts by A.~Arsie and P.~Lorenzoni.
\end{Remark}

Conversely, starting with a solution of (\ref{eq:genWDVV})--(\ref{eq:vphomo}),
it is possible to construct a Saito structure.

\begin{Proposition} \label{prop:fromvptoSaito}
Take constants $w_j\in\mathbb{C},$ $1\leq j\leq n$ satisfying $w_i-w_j\not\in\mathbb{Z}$ and $w_n=1$
and assume that $\vec{g}=(g_1,\dots,g_n)$ is a holomorphic solution of \eqref{eq:genWDVV}--\eqref{eq:vphomo}
 on a simply-connected domain $U$ in $\mathbb{C}^n$.
 Then there is a Saito structure on $U$ which has $(t_1,\dots,t_n)$ as a flat coordinate system.
 In addition, the Saito structure is semisimple $($i.e., $-\Phi(E)$ is semisimple$)$
 if and only if

$(SS)$ the $(n\times n)$-matrix $\big({-}(1+w_j-w_i) \frac{\partial g_j}{\partial t_i}\big)_{1\leq i,j\leq n}$
 is semisimple.
\end{Proposition}

\begin{proof}
 Define $E:=\sum_{i=1}^nw_it_i\partial_{t_i},$ $e:=\partial_{t_n},$
 $\mathcal{C}_{ij}:=\frac{\partial g_j}{\partial t_i},$
 $\tilde{\Phi}:={\rm d}\mathcal{C}$
 and $\knabla(\partial_{t_i})=0,$ $i=1,\dots,n$.
Then $(\knabla, \Phi, E, e)$ satisfies the conditions (a), (b) in Definition \ref{def:Saito}
and $\vec{g}$ becomes its potential vector field.
The last part of the proposition is obvious from $\mathcal{T}=-E\mathcal{C}$.
\end{proof}

Let $J$ be an $(n\times n)$-matrix with $J_{ij}=\delta _{i+j,n+1}$, $1\leq i,j \leq n$, where $\delta_{ij}$ denotes Kronecker's delta,
and, for an $(n\times n)$-matrix $A$, define $A^*$ by $A^*=JA^tJ$.

\begin{Proposition} \label{prop:conditionsFrob}
 Given a Saito structure on $M$,
 the following conditions are mutually equivalent:
 \begin{enumerate}\itemsep=0pt
 \item[$(i)$] For appropriate normalization of the flat coordinate system, it holds that $\mathcal{C}^*=\mathcal{C}$.
 \item[$(ii)$] For appropriate normalization of the flat coordinate system,
 there is a holomorphic function $F\in\mathcal{O}_M(U)$ such that
 \begin{equation*}
 \frac{\partial F}{\partial t_i}=g_{n+1-i}=(\vec{g}J)_i, \qquad i=1,\dots,n.
 \end{equation*}
 \item[$(iii)$] There is $r\in \mathbb{C}$ such that
 \begin{equation*}
 w_{n+1-i}+w_i=-2r,\qquad i=1,\dots,n,
 \end{equation*}
 and there is a metric $\eta$ $($in this paper, ``metric'' means non-degenerate symmetric $\mathbb{C}$-bilinear form on $TM)$
 such that
 \begin{gather*}
 \eta(X\star Y,Z)=\eta(X,Y\star Z), \\
 (\knabla\eta)(X,Y):={\rm d}(\eta(X,Y))-\eta(\knabla X,Y)-\eta(X,\knabla Y)=0, \\
 (E\eta)(X,Y):=E(\eta(X,Y))-\eta(EX,Y)-\eta(X,EY)=-2r\eta(X,Y),
 \end{gather*}
 for any $X,Y,Z\in \Theta_M$.
 \end{enumerate}
 A flat coordinate system has the ambiguity mentioned in Remark~{\rm \ref{rem:flat pvf constant}}.
 ``Normalization'' in the above conditions means to fix this ambiguity.
\end{Proposition}

\begin{proof} $({\rm ii})\Leftrightarrow ({\rm iii})$ is proved in \cite{Du, Sab}.
 $({\rm i})\Leftrightarrow ({\rm ii})$ is proved as follows.
 By definition, $\mathcal{C}^*=\mathcal{C}$ is equivalent to that $\mathcal{C}J$ is a symmetric matrix.
 From
 \begin{equation*}
 \frac{\partial g_{n+1-j}}{\partial t_i}=(\mathcal{C}J)_{ij},
 \end{equation*}
 we find that the symmetry of $\mathcal{C}J$ is equivalent to
 the existence of $F\in\mathcal{O}_M(U)$ such that $\frac{\partial F}{\partial t_i}=g_{n+1-i}$.
\end{proof}

The function $F$ in Proposition~\ref{prop:conditionsFrob} is called a {\it prepotential} in~\cite{Du} and a {\it potential} in~\cite{Sab}. It is well known that the prepotential satisfies the WDVV equation (cf.~\cite{Du}).

\subsection[Saito structure on the space of variables of isomonodromic deformations of an Okubo system]{Saito structure on the space of variables of isomonodromic deformations\\ of an Okubo system} \label{subsec:Okubo}

We consider an extended Okubo system (\ref{eq:okubopfaff}) with the same assumptions as in Section~\ref{sec:okuboiso}.
We show that the space of independent variables of (\ref{eq:okubopfaff}) can be equipped with
a Saito structure under some generic condition.
The arguments below closely follow \cite[Chapter VII]{Sab}.

In general setting, for a Saito bundle $(V,\nabla^V,\Phi^V,R_0^V,R_{\infty }^V)$ on a complex manifold $M$,
a $\nabla^V$-horizontal section $\omega$ of $V$ is said to be a {\it primitive section} if it satisfies the following conditions:

\begin{enumerate}\itemsep=0pt
 \item[(i)] $R_{\infty}^V\omega=\lambda \cdot\omega$ for some $\lambda\in\mathbb{C}$,
 \item[(ii)] $\varphi_{\omega}\colon TM \to V$ defined by $\varphi_{\omega}(X):=\Phi_X^V(\omega)$
 for $X\in\Theta_M$ is an isomorphism.
\end{enumerate}

Thanks to \cite[Chapter VII, Theorem 3.6]{Sab}, if there is a primitive section $\omega$,
we can introduce a Saito structure $(\knabla,\Phi,e,E)$ on $M$ via $\varphi_{\omega}$:
\begin{gather*}
 \knabla:=\varphi_{\omega}^{-1}\circ \nabla^V\circ \varphi_{\omega},\qquad
 \Phi:=\varphi_{\omega}^{-1}\circ \Phi^V\circ \varphi_{\omega},\qquad
 e:=\varphi_{\omega}^{-1}(\omega), \qquad E:=-\varphi_{\omega}^{-1}\big(R_0^V\omega\big).
\end{gather*}
In particular, we find
\[
 -\Phi(E)=\varphi_{\omega}^{-1}\circ R_0^V\circ \varphi_{\omega}, \qquad
 \nabla E=\varphi_{\omega}^{-1}\circ R_{\infty}^V\circ \varphi_{\omega}^{-1}+(1-\lambda){\rm id}_V,
\]
where $\lambda$ is the eigenvalue of $\omega$ in the condition (i).

Returning to our setting,
we consider an extended Okubo system (\ref{eq:okubopfaff}) defined on $\mathbb{P}^1\times U$,
where $U$ is a domain of $\mathbb{C}^n$.
For a while, we treat (\ref{eq:okubopfaff}) on an appropriate smaller domain $W\subset U\setminus\Delta_H$
 so that we can take the coordinate system $(z_1,\dots,z_n)$ and an invertible matrix $P$ such that
\begin{equation} \label{eq:PtoT,B}
 P^{-1}TP=\operatorname{diag}[ z_1,\dots, z_n],\qquad
 P^{-1}\tilde{\Omega}P
 =-\operatorname{diag} [{\rm d}z_1,\dots, {\rm d}z_n ].
\end{equation}
Then, as we observed in Section~\ref{sec:okuboiso}, the extended Okubo system (\ref{eq:okubopfaff}) induces
a Saito bundle $\big(W\times \mathbb{C}^n, {\rm d}, \alpha\tilde{\Omega}^t,T^t, B_{\infty}\big)$
with some $\alpha\in\mathbb{C}\setminus\{0\}$.
Any constant eigenvector of the matrix $B_{\infty}$ is ${\rm d}$-horizontal
and satisfies the condition~(i).
We take ${\bf e}_n$ as such an eigenvector without loss of generality.
(We can reduce other cases to this one applying a gauge transformation to the extended Okubo system~(\ref{eq:okubopfaff})
by a constant matrix.)

\begin{Lemma} \label{lem:nscondflat}
 The constant eigenvector ${\bf e}_n$ of $B_{\infty}$ is a primitive section
 of the Saito bundle $\big(W\times\mathbb{C}^n, {\rm d}, \alpha\tilde{\Omega}^t, \alpha T^t, B_{\infty}\big)$
 if and only if $P_{nj}\neq 0$, $1\leq j\leq n$, at any point on $W$.
\end{Lemma}

\begin{proof}
 For the Saito bundle $\big(W\times\mathbb{C}^n, {\rm d}, \alpha\tilde{\Omega}^t, \alpha T^t, B_{\infty}\big)$,
 the definition of the bundle map $\varphi_{{\bf e}_n}$ reads $\varphi_{{\bf e}_n}(X)=\alpha\tilde{\Omega}_X^t{\bf e}_n$.
 Noting ${\bf e}_n=({\bf e}_1,\dots,{\bf e}_n)\big(P^t\big)^{-1}(P_{n1},\dots,P_{nn})^t$ and~(\ref{eq:PtoT,B}), we have
 \[
 \varphi_{{\bf e}_n}(X)=-\alpha({\bf e}_1,\dots,{\bf e}_n)\big(P^t\big)^{-1}\operatorname{diag}[P_{n1},\dots,P_{nn}](Xz_1,\dots,Xz_n)^t,
 \qquad X\in\Theta_M.
 \]
 Then it is clear that $\varphi_{{\bf e}_n}$
 is an isomorphism if and only if $\prod_{j=1}^nP_{nj}\neq 0$.
\end{proof}

\begin{Remark}
 We consider the Saito bundle $\big(W\times\mathbb{C}^n, {\rm d}, \alpha\tilde{\Omega}^t, \alpha T^t, B_{\infty}\big)$.
 Let us assume~(\ref{cond:T}).
 Then we have $e=\varphi_{{\bf e}_n}^{-1}({\bf e}_n)=\alpha^{-1}\partial_{x_n}$.
 So, in this case, it is natural to choose $\alpha=1$.
\end{Remark}

In the following, we always choose $\alpha=1$.

\begin{Lemma} \label{lem:PtoT}
 For an extended Okubo system~\eqref{eq:okubopfaff},
 the following two conditions are equivalent:
 \begin{enumerate}\itemsep=0pt
 \item[$(i)$] ${\rm d}T_{n1}\wedge \cdots \wedge {\rm d}T_{nn}\neq 0$
 at any point on $W$,
 \item[$(ii)$] $\prod_{j=1}^nP_{nj}\neq 0$ at any point on~$W$.
 \end{enumerate}
\end{Lemma}

\begin{proof} From (\ref{eq:ci1}), it holds that
 \[
 {\rm d}T_{nj}=(\lambda_n-\lambda_j-1)\tilde{\Omega}_{nj},
 \]
 from which and (\ref{eq:taikakuBi}) we have
 \begin{gather*}
 ({\rm d}T_{n1},\dots,{\rm d}T_{nn}) \operatorname{diag}[\lambda_n-\lambda_1-1,\lambda_n-\lambda_2-1,\dots,-1]^{-1}P\\
 \qquad {} =\big(\tilde{\Omega}_{n1},\dots,\tilde{\Omega}_{nn}\big)P 
 =-({\rm d}z_1,\dots,{\rm d}z_n) \operatorname{diag}[P_{n1},\dots,P_{nn}].
 \end{gather*}
 Hence we obtain that (i)$\iff$(ii).
\end{proof}

So far we have treated (\ref{eq:okubopfaff}) on $W$ outside of $\Delta_H$.
From now on, we consider (\ref{eq:okubopfaff}) on $U$ including $\Delta_H$.

\begin{Theorem} \label{saitojacobian}
 An extended Okubo system \eqref{eq:okubopfaff} induces
a Saito structure on~$U$ if and only if
 \begin{equation} \label{jacobian}
 {\rm d}T_{n1}\wedge \cdots \wedge {\rm d}T_{nn} \neq 0 \qquad \text{on} \ U.
 \end{equation}
 If \eqref{jacobian} is satisfied,
 the set of variables $t_j:=-(\lambda_j-\lambda_n+1)^{-1}T_{nj}=C_{nj}$, $1\leq j\leq n$ gives a~flat coordinate system
 on $U$.
\end{Theorem}

\begin{proof}
 In virtue of Lemmas~\ref{lem:nscondflat} and~\ref{lem:PtoT},
 the theorem holds on $W\subset U\setminus \Delta_H$.
 We also see that $\{t_j=-(\lambda_j-\lambda_n+1)^{-1}T_{nj}\}$ is a flat coordinate system by Lemma~\ref{lem:flatC=t}.
 Then $(\knabla,\Phi,E,e)$, where
 $\Phi:={\rm d}C$, $E:=\sum_{k=1}^n(\lambda_k-\lambda_n+1)t_k\partial_{t_k}$, $e:=\partial_{t_n}$ and~$\knabla$ is a connection defined by $\knabla(\partial_{t_i})=0$,
 satisfies the conditions~(a),~(b) of Saito structure on~$W$.
 Due to the identity theorem,
 $(\knabla,\Phi,E,e)$ satisfies the conditions (a),~(b) on $U$.
 Hence~(\ref{eq:okubopfaff}) induces a Saito structure on $U$.
\end{proof}

\begin{Remark} \label{rem:MSaitoSabbah}
 Theorem~\ref{saitojacobian} is close to a result of Morihiko Saito in~\cite{SaM},
 which (together with a mixed Hodge structure) implies existence of Kyoji Saito's primitive forms
 in the case of non-quasihomogeneous isolated hypersurface singularities.
 The book of Sabbah \cite{Sab} contains in Chapter~VII an improved and easier accessible
 (than the result in \cite{SaM}) version of M.~Saito's result.
\end{Remark}

\begin{Remark} \label{rem:uniqunessofSaitoOkubo}
 In general, the Saito structure induced by an extended Okubo system depends on the choice of constant eigenvectors of $B_{\infty}$.
 However, in the case where all of the eigenvalues $\{\lambda_1,\dots,\lambda_n\}$ are real numbers
 and satisfy $\lambda_1\leq\cdots\leq\lambda_{n-1}<\lambda_n$,
 we can distinguish $\lambda_n$ and the eigenvector belonging to it from the remaining eigenvectors,
 and thus a unique Saito structure can be specified.
 This is the origin of the uniqueness of the Saito structure on the orbit space of a well-generated unitary reflection group
 in Section~\ref{sec:complexreflection}.
\end{Remark}

 We consider the reduced form (\ref{eq:okubopfaff0}) in Remark~\ref{Re:z=xn} of an extended Okubo system
 with $\det B_{\infty}\neq 0$ taking a flat coordinate system $(t_1,\dots,t_n)$:
 \begin{equation} \label{eq:reduced2}
 {\rm d}Y_0=\left(\sum_{i=1}^n T^{-1}\tilde{B}^{(i)}B_{\infty}\,{\rm d}t_i\right)Y_0.
 \end{equation}

\begin{Corollary} \label{cor:tandokuJ}
 Let $V_i$, $1\leq i\leq n$, be the vector fields defined by \eqref{eq:Saito logvec}.
 Then it holds that
 \begin{equation*} 
 Y_0=\big(y_1^{(0)},\dots,y_n^{(0)}\big)^t=-B_{\infty}^{-1}\big(V_ny_n^{(0)},\dots,V_1y_n^{(0)}\big)^t.
 \end{equation*}
 This means the ``primitivity'' of $y_n^{(0)}$ in the sense of K.~Saito.
\end{Corollary}

\begin{proof}
 From (\ref{eq:reduced2}), (\ref{eq:commute}) and (\ref{eq:ci1}), we have
 \begin{equation*}
 \frac{\partial y_n^{(0)}}{\partial t_i}=\big(\tilde{B}^{(i)}_{n1},\dots,\tilde{B}^{(i)}_{nn}\big)T^{-1}B_{\infty}Y_0
 =\big(\tilde{B}^{(n)}_{i1},\dots,\tilde{B}^{(n)}_{in}\big)T^{-1}B_{\infty}Y_0.
 \end{equation*}
 Noting $\tilde{B}^{(n)}_{ij}=\delta_{ij}$, we obtain
 \begin{equation*}
 \big(V_ny_n^{(0)},\dots,V_1y_n^{(0)}\big)^t=-T\big(\partial_{t_1}y^{(0)}_n,\dots,\partial_{t_n}y_n^{(0)}\big)^t
 =-TT^{-1}B_{\infty}Y_0=-B_{\infty}Y_0.\tag*{\qed}
 \end{equation*}\renewcommand{\qed}{}
\end{proof}

\begin{Corollary} \label{cor:corgexWDVVandPVI}
 We consider the case of $n=3$.
 There is a correspondence between generic solutions satisfying the semisimplicity condition $(SS)$ in Proposition~{\rm \ref{prop:fromvptoSaito}} to the extended WDVV equation
 \begin{gather*}
 \sum_{m=1}^3\frac{\partial^2 g_m}{\partial t_k\partial t_i}\frac{\partial^2 g_j}{\partial t_l\partial t_m}=
 \sum_{m=1}^3\frac{\partial^2 g_m}{\partial t_l\partial t_i}\frac{\partial^2 g_j}{\partial t_k\partial t_m},\qquad i,j,k,l=1,2,3, \\
 \frac{\partial ^2 g_j}{\partial t_3\partial t_i}=\delta_{ij},\qquad i,j=1,2,3, \\
 Eg_j=\sum_{k=1}^3w_kt_k\frac{\partial g_j}{\partial t_k}=(1+w_j)g_j,\qquad j=1,2,3.
 \end{gather*}
 and generic solutions to the Painlev\'e~VI equation
 \begin{gather*}
 \frac{{\rm d}^2y}{{\rm d}t^2}= \frac{1}{2}\left(\frac{1}{y}+\frac{1}{y-1}+\frac{1}{y-t}\right)\left(\frac{{\rm d}y}{{\rm d}t}\right)^2
 -\left(\frac{1}{t}+\frac{1}{t-1}+\frac{1}{y-t}\right)\frac{{\rm d}y}{{\rm d}t} \\
\hphantom{\frac{{\rm d}^2y}{{\rm d}t^2}=}{} +\frac{y(y-1)(y-t)}{t^2(t-1)^2}\left(\alpha+\beta\frac{t}{y^2}+\gamma\frac{t-1}{(y-1)^2}+\delta\frac{t(t-1)}{(y-t)^2}\right).
 \end{gather*}
 Let $T$ be the $(3\times 3)$-matrix whose entries are given by
 \begin{equation*}
 T_{ij}=-(1-w_i+w_j)\frac{\partial g_j}{\partial t_i},
 \end{equation*}
 and take $P$ such that $P^{-1}TP$ is a diagonal matrix.
 Let
 \begin{equation*}
 r_i=-\big(P \operatorname{diag}[w_1-w_3, w_2-w_3, 0]P^{-1}\big)_{ii},\qquad i=1,2,3
 \end{equation*}
 and $\theta_{\infty}=w_1-w_2$.
 Then the correspondence between the parameters is given by
 \begin{equation*}
 \alpha=\frac{1}{2}(\theta_{\infty}-1)^2,\qquad \beta=-\frac{1}{2}r_1^2,\qquad \gamma=\frac{1}{2}r_2^2,\qquad \delta=\frac{1}{2}\big(1-r_3^2\big).
 \end{equation*}
\end{Corollary}

\begin{proof}
 Since the condition (\ref{jacobian}) in Theorem \ref{saitojacobian} is generic one,
 we have a correspondence between generic extended Okubo systems and solutions to the extended WDVV equation
 satisfying the semisimplicity condition.
 Extended Okubo systems can be regarded as isomonodromic deformations of Okubo systems
 and, in the rank-three case, it is known that such isomonodromic deformations are governed by generic solutions to the Painlev\'e VI equation.
 Indeed, the reduction of a $3\times 3$ Okubo system to a $2\times 2$ system with four regular singular points
 on $\mathbb{P}^1$ and the reconstruction of a $3\times 3$ Okubo system from a $2\times 2$ system
 with four regular singular points are explained in \cite{BoC,Bo0,Bosix}.
\end{proof}

\begin{Remark}
 In the papers \cite{AL0,Lo},
 Arsie and Lorenzoni treated the relationship between semisimple bi-flat $F$-manifolds and the sigma form of Painlev\'e VI equation
 via generalized Darboux--Egorov systems.
In particular, \cite{Lo} shows that the three-dimensional regular semisimple bi-flat $F$-manifolds
 are parameterized by solutions to the (full-parameter) Painlev\'e VI equation.
Recently it is proved in~\cite{AL2,KoMiSh} that bi-flat $F$-manifold is equivalent to Saito structure.
Therefore Corollary~\ref{cor:corgexWDVVandPVI} provides another proof of Arsie--Lorenzoni's result.
The proof here makes clear the relationship between Saito structures and isomonodromic deformations of an Okubo system.

 In \cite{AL1}, Arsie and Lorenzoni study the relationship between three-dimensional regular non-semisimple bi-flat $F$-manifolds
 and Painlev\'e~IV and~V equations.
 Then it is naturally expected that there is a correspondence between solutions to the extended WDVV equation
not satisfying the semisimple condition~(SS)
 and isomonodromic deformations of generalized Okubo systems introduced in~\cite{Kaw}.
 This problem is studied in~\cite{KaMa}.
\end{Remark}

\section{Flat basic invariants for complex reflection group} \label{sec:complexreflection}

Let $G$ be an irreducible well-generated complex reflection group.
The main purpose of this section is to prove the existence of distinguished basic invariants
of $G$ called {\it flat basic invariants} (Theorem~\ref{prop:C_1} and Definition~\ref{def:flatbasicinv}).
For this purpose,
we construct an extended Okubo system on the orbit space of $G$
 which is called the $G$-{\it quotient system}.
 The flat basic invariants of $G$ are defined by the flat coordinate system of the Saito structure on the orbit space
induced by the $G$-quotient system
applying Theorem~\ref{saitojacobian}.
This is a natural extension of K. Saito's flat basic invariants for real reflection groups~\cite{Sai2,SYS}.
As a consequence of Theorem~\ref{prop:C_1}, we find that the potential vector field
for any well-generated complex reflection group $G$ has polynomial entries
(Corollary~\ref{cor:polyG}).
It is underlined that our proof of Theorem~\ref{prop:C_1} is constructive:
we obtain a~method to explicitly construct the flat basic invariants and the potential vector field for~$G$.
We explain it for an example (Example~\ref{ex:G26}).
Explicit representations of potential vector fields for exceptional complex reflection groups are found in~\cite{KMS}.
(There is a procedure to construct flat basic invariants from the potential vector fields
for well-generated complex reflection groups, see Remark~\ref{rem:KMS}.)

Let $G$ be a finite irreducible complex (unitary) reflection group
acting on the standard representation space
\[
U_n=\{u=(u_1,u_2,\dots, u_{n})\,|\, u_j\in\mathbb{C} \},
\]
and let $F_{i}(u)$ of degree $d_i$, $1\le i\le n$ be a fundamental system of $G$-invariant homogeneous polynomials.
We assume that
\[
 d_1\le d_2\le \dots \le d_{n}.
\]
We define coordinate functions on $X:=U_n/G$ by,
$x_i=F_{i}(u)$, $1\le i\le n$.
Let $D\subset X$ be the branch locus of $\pi_G\colon U_n\to X$.
Let $h(x)$ be the (reduced) defining function of $D$ in the coordinate $x=(x_1,x_2,\dots,x_n)$.
We assume
that $G$ is well generated (see, e.g.,~\cite{Be}).
Then it is known that~$h(x)$ is a monic
polynomial in~$x_n$ of degree~$n$ (cf.\ Remark~\ref{rem:flatvector}).
We define a~weight~$w(\cdot)$ by~$w(x_i)=d_i/d_n$.
Then~$h(x)$ is a weighted homogeneous polynomial in~$x$.
It is known that~$D$ is free~\cite{OT,Tr1}.
Here we give a fundamental lemma for the discriminant~$h(x)$.

\begin{Lemma}\label{lemma:B4 of h} $H(z):=h(x',x_n+z)$ satisfies the assumption $(\rm{A}3)$.
\end{Lemma}

\begin{proof}It is known \cite{Be} that there is a system of generators $\{ V_1,\dots ,V_n\}$ of the free $\mathbb{C}[x]$-module $\operatorname{Der}(-\log D)$
such that $V_{i}=\sum_{j=1}^n v_{ij}(x) \,\partial_{x_j}$, $1\le i\le n$ are weighted homogeneous and satisfy
$[\partial_{x_n},V_i-x_n\partial_{x_n+1-i}]=0.$
In particular, it holds that
\[
 \deg_{x_n}(v_{ij}(x))=
\begin{cases}
 1\quad \mbox{if}\ i+j=n+1,\\
 0\quad \mbox{if}\ i+j\ne n+1.
\end{cases}
\]

In virtue of Lemma~\ref{lemma:zeroV}, it is sufficient to show that $V=0$
for any $V=\sum_{i=1}^{n} c_i(x')\,\partial_{x_i}\in\operatorname{Der}(-\log D)$.
There are weighted homogeneous polynomials $a_i(x)\in \mathbb{C}[x]$ such that
$V=\sum_{i=1}^n a_i(x) V_i$, which is equivalent to
\begin{equation}\label{equ:c_j}
c_j(x')=\sum_{i=1}^n a_i(x) v_{ij}(x), \qquad 1\le j\le n.
\end{equation}
Let $I=\{ i\,|\,1\leq i\leq n,\, a_i(x)\ne 0\}$, and assume $I\ne \varnothing$.
Let $i_0\in I$ be so taken that
\[
 \deg_{x_n}(a_{i_0})\ge \deg_{x_n}(a_{i}),\qquad \forall\, i\in I.
\]
Then
\begin{gather*}
 \deg_{x_n}(a_{i_0}(x) v_{i_0, n+1-i_0}(x))\\
 \qquad{} =\deg_{x_n}(a_{i_0}(x))+1
 >\deg_{x_n}(a_{i}(x) v_{i,n+1-i_0}(x)),\qquad i\ne i_0,\quad i\in I,
\end{gather*}
which implies
\begin{gather*}
 \deg_{x_n}\left(\sum_{i=1}^n a_i(x) v_{i,n+1-i_0}(x) \right)
 =\deg_{x_n}(a_{i_0}(x) v_{i_0, n+1-i_0}(x))>0
 =\deg_{x_n} (c_{n+1-i_0}(x')).
\end{gather*}
This contradicts (\ref{equ:c_j}) with $j=n+1-i_0$,
which implies that $I=\varnothing$, that is, $V=0$.
\end{proof}

 Let $H_i:=\{\ell_i(u)=0\}\subset U_n$, $1\le i\le N$ be all the
 distinct reflecting hyperplanes of $G$.
 If~$H_i$ is the reflecting hyperplane of a unitary reflection
 $g_i\in G$ of the order $m_i$, then we have
 $g_i(H_i)=H_i$, $g_j(H_i)\ne H_i$ for $j\ne i$. As stated in \cite{ST}
 \begin{equation}\label{equ:h(x(u))}
 h(x(u))=c_1 \prod_{i=1}^N \ell_i(u)^{m_i},\qquad
 \det \left({\partial x\over \partial u}\right)
 =c_2 \prod_{i=1}^N \ell_i(u)^{m_i-1},
 \end{equation}
 for some $c_1,c_2\in \mathbb{C}^{\times}$,
 where $ \big({\partial x\over \partial u}\big)
 =\big({\partial x_j\over \partial u_i}\big)_{i,j=1,2,\dots,n}$.

Put
\[
 \tilde{h}(u)=\prod_{i=1}^N \ell_i(u),\qquad
 \text{and}\qquad \tilde{D}=\cup_{i=1}^N H_i.
\]

Let
\[
 M_{\mathcal{V}}(x)=M_{\mathcal{V}^{(h)}}(x),\qquad
 V_i(x)=V_i^{(h)}(x),\qquad 1\le i\le n,
\]
where $M_{\mathcal{V}^{(h)}}(x)$, $V_i^{(h)}(x)$ are respectively defined by
 (\ref{equ:Def of M_V and M_h}) and (\ref{equ:LogVFs}) with respect to~$h(x)$.
Recall that $\det(M_{\mathcal{V}}(x))=h(x)$, and
$V_1(x)$ is the Euler operator:
$V_1(x)=\sum_{i=1}^n w(x_i) x_i \partial_{x_i}$.

 Let $\tilde{V}_i$ be the pull-back of $V_i$ to $U_n$, that is,
 \begin{gather} \label{eq:basicdiffU}
 \begin{pmatrix}
 \tilde{V}_n&\tilde{V}_{n-1}& \dots& \tilde{V}_1
 \end{pmatrix}^t
 =M_{\tilde{\mathcal{V}}}
 \begin{pmatrix} \partial_{u_1}& \partial_{u_2}& \dots&
 \partial_{u_n} \end{pmatrix}^t,\qquad\! \text{where}\quad\!
 M_{\tilde{\mathcal{V}}}=M_{\mathcal{V}}\left({\partial x\over \partial u}\right)^{-1}\!.\!\!\!
 \end{gather}
Then, from (\ref{equ:h(x(u))}), it holds that
 \begin{equation*}
 \det(M_{\tilde{\mathcal{V}}})=c \tilde{h}(u),
 \end{equation*}
for a constant number $c\ne 0$.

In what follows, we prove the existence of flat basic invariants for the well-generated complex reflection group $G$
along the following line:
\begin{description}\itemsep=0pt
 \item[\rm Step 1.]
 For any homogeneous linear function $y(u)$ of $u$, we show that the Pfaffian system (\ref{equ:G-quotient})
 satisfied by $\hat{Y}=\big(\tilde{V}_ny,\dots, \tilde{V}_1y\big)^t$ descends on $X=U_n/G$ (in Lemma~\ref{lemma:L1}).

 \item[\rm Step 2.]
 We show that the Pfaffian system (\ref{equ:G-quotient}) is transformed into the reduced form
 of an extended Okubo system (\ref{equ:w(T)=w(Mh)})
 applying a gauge transformation by an upper triangular matrix (in Lemma~\ref{lemma:w(T)=w(Mh)}).
 The resulting extended Okubo system (\ref{equ:w(T)=w(Mh)}) is the $G$-quotient system, see Remark~\ref{rem:gs vs okubo}.

 \item[\rm Step 3.]
 Applying Theorem~\ref{saitojacobian} to the $G$-quotient system (\ref{equ:w(T)=w(Mh)}),
 we obtain a flat coordinate system $t=(t_1,\dots,t_n)$ on $X=U_n/G$
 and show that $t$ regarded as a system of functions on $U_n$ gives a system of generators of $G$-invariants,
 which is nothing but the flat basic invariants of $G$ (in Theorem~\ref{prop:C_1}).
\end{description}

We start with showing a lemma which is used in the proof of Lemma~\ref{lemma:L1}:

\begin{Lemma}\label{lemma:tildeVi is log}\quad
\begin{enumerate}\itemsep=0pt
 \item[$(i)$]
 The entries of $M_{\tilde{\mathcal{V}}}$ are polynomials in $u$,
 $\tilde{V}_k$ $(1\le k\le n)$
 are logarithmic vector fields along~$\tilde{D}$, and $\tilde{D}$ is free.

\item[$(ii)$]
 For any $k\in \{ 1,2,\dots,n\}$, all entries of
 $\big(\tilde{V}_k M_{\tilde{\mathcal{V}}}\big) M_{\tilde{\mathcal{V}}}^{-1}$
 are $G$-invariant polynomials in $u$, that is,
 polynomials in $x$.
 \end{enumerate}
\end{Lemma}

\begin{proof}
 (i) is known by Terao and others \cite{OT,Tr1}.
 We prove (ii).
 Fix $i\in \{1,\dots, N\}$ arbitrarily, and we shall prove
 $\big(\tilde{V}_k M_{\tilde{\mathcal{V}}}\big)M_{\tilde{\mathcal{V}}}^{-1}$ is holomorphic along $H_i$.
 Let $y_1(u)=\ell_i(u)$, and choose $y_k(u)=\ell_{i_k}(u)$,
 $2\le k\le n$,
 so that $y_1,y_2,\dots,y_n$ are linearly independent.
 Then (i) implies that
 \begin{equation*}
 \tilde{V}_k=\sum_{j=1}^n v_{k,j}(y) y_j\partial_{y_j},\qquad 1\le k\le n,
 \end{equation*}
 for some $v_{k,j}(y)\in \mathbb{C}[y]$.
 This is equivalent to
 \begin{equation}\label{equ: V_k h 4}
 M_{\tilde{\mathcal{V}}}=M'_{\tilde{\mathcal{V}}}(y)\cdot
 \operatorname{diag}[y_1,y_2,\dots,y_n],
 \end{equation}
 for some $M'_{\tilde{\mathcal{V}}}(y)\in \mathbb{C}[y]^{n\times n}$.
 Then it also holds that $\tilde{V}_k\,M_{\tilde{\mathcal{V}}}=M''_{\tilde{\mathcal{V}}}(y)\cdot
 \operatorname{diag}[y_1,y_2,\dots,y_n]$,
 for some $M''_{\tilde{\mathcal{V}}}(y)\in \mathbb{C}[y]^{n\times n}$.
 Thus we have
 \begin{equation}\label{equ: V_k h 5}
 \big(\tilde{V}_k M_{\tilde{\mathcal{V}}}\big)M_{\tilde{\mathcal{V}}}^{-1}
 =M''_{\tilde{\mathcal{V}}}(y) \big(M'_{\tilde{\mathcal{V}}}(y)\big)^{-1}.
 \end{equation}
 Since $\det(M_{\tilde{\mathcal{V}}})=\prod_{j=1}^N \ell_j(u)$, the equality~(\ref{equ: V_k h 4})
 implies that
 $\big(M'_{\tilde{\mathcal{V}}}(y)\big)^{-1}$
 is holomorphic along $\{ y_1=0\}$.
 Consequently the equality~(\ref{equ: V_k h 5})
 implies that $\big(\tilde{V}_k M_{\tilde{\mathcal{V}}}\big)M_{\tilde{\mathcal{V}}}^{-1}$
 is holomorphic along $H_i$.
 This proves that $\big(\tilde{V}_k M_{\tilde{\mathcal{V}}}\big) M_{\tilde{\mathcal{V}}}^{-1}\in
 \mathbb{C}[u]^{n\times n}$.
 It is easy to see that $\big(\tilde{V}_k\,M_{\tilde{\mathcal{V}}}\big) M_{\tilde{\mathcal{V}}}^{-1}$ is $G$-invariant.
\end{proof}

Put
 \begin{equation} \label{eq:hat Y}
 \hat{Y}=\begin{pmatrix}\hat{y}_1&\hat{y}_2& \dots& \hat{y}_n \end{pmatrix}^t
 :=\begin{pmatrix} \tilde{V}_n&\tilde{V}_{n-1}& \dots& \tilde{V}_1 \end{pmatrix}^t y(u),
 \end{equation}
 for any homogeneous linear function $y(u)$ of $u$.
Then we have
\begin{equation*}
 w(\hat{y}_i)=w(V_{n-i+1})+(1/d_n)=(d_n-d_i+1)/d_n,\qquad 1\le i\le n.
 \end{equation*}
 Let
 \[
 \big({\rm d} M_{\tilde{\mathcal{V}}}\big)M_{\tilde{\mathcal{V}}}^{-1}
 =\sum_{k=1}^n\hat{B}^{(k)}\,{\rm d}x_k.
 \]
 Then by (\ref{eq:basicdiffU}), it is clear that (\ref{eq:hat Y}) satisfies the Pfaffian system
 \begin{equation}\label{equ:G-quotient}
{\rm d}\hat{Y}=\left[\sum_{k=1}^n\hat{B}^{(k)}\, {\rm d}x_k\right]\hat{Y}.
\end{equation}

\begin{Lemma}\label{lemma:L1}
 All entries of $h(x)\hat{B}^{(k)}$ are weighted homogeneous polynomials
 in $x$, that is,
\eqref{equ:G-quotient}~is defined on $X=U_n/G$ with its singular locus along~$D$.
\end{Lemma}

\begin{proof}
 Put $P_k(x)=\big(\tilde{V}_k M_{\tilde{\mathcal{V}}}\big) M_{\tilde{\mathcal{V}}}^{-1}$ for
 $k=1,2,\dots n$. Then $P_k(x)\in \mathbb{C}[x]^{n\times n}$ from (ii) of
 Lemma~\ref{lemma:tildeVi is log}.
 Using $\big(\tilde{V}_n,\dots,\tilde{V}_1\big)^t=M_{\mathcal{V}}(\partial_{x_1},\dots,\partial_{x_n})^t$,
 it is clear that
 \[ \hat{B}^{(k)}=\sum_{j=1}^n\big(M_{\mathcal{V}}^{-1}\big)_{k,j} P_{n-j+1}(x),\]
 which implies
 $h(x)\hat{B}^{(k)}(x)\in \mathbb{C}[x]^{n\times n}$.
\end{proof}

Put
\begin{equation}\label{equ:B_inf}
 B_{\infty}:=\operatorname{diag}\left[w(x_1),w(x_2),\dots,w(x_n)\right]
-\left(1+{1\over d_n}\right)I_n.
\end{equation}

\begin{Lemma}\label{lemma:w(T)=w(Mh)} Let $\hat{Y}$, $\hat{B}^{(k)}$ be the same as in Lemma~{\rm \ref{lemma:L1}}.
 Then there is an upper triangular matrix $R(x')\in {\rm GL}(n,\mathbb{C}[x'])$ satisfying $R(0)=I_n$
 such that, if we put
 \begin{gather*}
 Y=(y_1,y_2,\dots,y_n)^t:=-B_{\infty}^{-1}R(x')\hat{Y}, \\
 B^{(k)}:=B_{\infty}^{-1}\left(R(x')\hat{B}^{(k)}R(x')^{-1}+{\partial R(x')\over \partial x_k}R(x')^{-1}\right)B_{\infty},\qquad
 1\le k\le n,
 \end{gather*}
 then the system of differential equations
 \begin{equation}\label{equ:w(T)=w(Mh)}
 {\rm d}Y=\left[\sum_{k=1}^nB^{(k)}\,{\rm d}x_k\right]Y
 \end{equation}
 satisfied by $Y$ is the reduced form of an extended Okubo system $($see Remark~{\rm \ref{Re:z=xn})},
 and the residue matrix $\big(B^{(n)}\big)_{\infty}$ of $B^{(n)}\,{\rm d}x_n$ at $x_n=\infty$
 is equal to the diagonal matrix~$B_{\infty}$ in~\eqref{equ:B_inf}.
 The matrix $R(x')$ may be chosen so that all the entries
 are weighted homogeneous polynomials in $x'$
 with $w(R(x')_{i,j})=w(\hat{y_i})-w(\hat{y_j})$.
\end{Lemma}

\begin{proof}
 Let $\hat{B}_{\infty}(x')=-\big(x_n\hat{B}^{(n)}\big)|_{x_n=\infty}$ be the residue
 matrix of $\hat{B}^{(n)}\,{\rm d}x_n$ at $x_n=\infty$.
 From
 \begin{gather*} 
 w\big(\hat{B}^{(k)}_{i,j}(x')\big) =w(\hat{y}_i)-w(\hat{y}_j)-w(x_k)\\
 \hphantom{w\big(\hat{B}^{(k)}_{i,j}(x')\big)}{} =w(V_{n-i+1})-w(V_{n-j+1})-w(x_k)
 =-w(x_i)+w(x_j)-w(x_k),
 \end{gather*}
 we find that the degree of
 $h(x)\hat{B}^{(n)}_{i,j}(x)$ in $x_n$ is at most $n-1$,
 which implies that $-\hat{B}_{\infty}(x')$ is the coefficients of
 $x_n^{n-1}$
 of $h(x)\hat{B}^{(n)}(x)$.
 Consequently $\big(\hat{B}_{\infty}\big)_{i,j}$ is a weighted homogeneous
 polynomial
 in $x'$ with
 \[
 w\big(\big(\hat{B}_{\infty}\big)_{i,j}\big)=w(x_j)-w(x_i).
 \]

 Let $n_1,n_2,\dots ,n_k$ be the positive integers such that
 $n_1+n_2+\dots +n_k=n$ and
 \[
 w(x_1)=\dots =w(x_{n_1})<w(x_{n_1+1})=\dots =w(x_{n_1+n_2})<\dots
 \le w(x_n).
 \]

Then
 $\hat{B}_{\infty}(x')$ has the form
 \[
 \hat{B}_{\infty}(x')= \begin{pmatrix}
 R_1&*&\cdots &*\\
 O&R_2&\cdots &*\\
 \cdots&\cdots&\cdots &\cdots\\
 O&O&\cdots &R_k
 \end{pmatrix},
 \]
 for some $R_i\in \mathbb{C}^{n_i\times n_i}$, $1\le i\le k$.

 Since
 \[
 \left(\sum_{k=1}^nw(x_k)x_k\hat{B}^{(k)}(x)\right)\hat{Y}
 =\left(\sum_{k=1}^nw(x_k)x_k{\partial\over \partial x_k}\right)\hat{Y}
 =
 \operatorname{diag}\big[w(\hat{y}_1),w(\hat{y}_2),\dots,w(\hat{y}_n)\big]
 \hat{Y},
 \]
 for all solutions $\hat{Y}$ of (\ref{equ:G-quotient}),
 we have
 \begin{equation} \label{eq:Eulerdiag}
 \sum_{k=1}^nw(x_k)x_k\hat{B}^{(k)}(x)=
 \operatorname{diag}\big[w(\hat{y}_1),w(\hat{y}_2),\dots,w(\hat{y}_n)\big].
 \end{equation}
 Substituting $x'=0$ for (\ref{eq:Eulerdiag}), we have
 \begin{equation*}
 \big(x_n\hat{B}^{(n)}(x)\big)|_{x'=0}=\operatorname{diag}\big[w(\hat{y}_1),w(\hat{y}_2),
 \dots,w(\hat{y}_n)\big],
 \end{equation*}
 which shows $R_i$ are diagonal, and $\hat{B}_{\infty}(x')$ is
 an upper triangular matrix with the diagonal elements $-w(\hat{y}_i)$.
 Then, by elementary linear algebra, we find that there is
 an upper triangular matrix $R(x')\in {\rm GL}(n,\mathbb{C}[x'])$
 with the form
 \begin{equation} \label{eq:R}
 R(x')= \begin{pmatrix}
 I_{n_1}&*&\cdots &*\\
 O&I_{n_2}&\cdots &*\\
 \cdots&\cdots&\cdots &\cdots\\
 O&O&\cdots &I_{n_k}
 \end{pmatrix},
 \end{equation}
 and satisfying
 \[
 R(x')\hat{B}_{\infty}(x')R(x')^{-1}=
 -\operatorname{diag}\big[w(\hat{y}_1),w(\hat{y}_2),\dots,
 w(\hat{y}_n)\big]
 =B_{\infty}.
 \]
 By construction of $R(x')$, we find that all entries $R(x')_{i,j}$ are
 weighted homogeneous with $w(R(x')_{i,j})=w(\hat{y_i})-w(\hat{y_j})$.
 Now put $B^{(k)}=B_{\infty}^{-1}\big(R(x')\hat{B}^{(k)}R(x')^{-1}
 +{\partial R(x')\over \partial x_k}R(x')^{-1}\big)B_{\infty},$ $1\le k\le n$
 and $B=\sum_{k=1}^nB^{(k)}\,{\rm d}x_k$.
 Then we find that
 \begin{equation*}
 \big(B^{(n)}\big)_{\infty}=-\big(x_n B^{(n)}\big)|_{x_n=\infty}
 =B_{\infty}^{-1}R(x')\hat{B}_{\infty}R(x')^{-1}B_{\infty}=B_{\infty}.
 \end{equation*}
 Here we note that the inverse image of $h(x',x_n)=\prod_{i=1}^n\big(x_n-z_i^0(x')\big)=0$ by the quotient mapping~$\pi_G$
 is the union of all reflecting hyperplanes for $G$.
 Then, at a generic point~$x'$, we see, from~(\ref{eq:hat Y}),
 that the local monodromy matrix of~(\ref{equ:G-quotient}) around $x_n=z_i^0(x')$ $(i=1,\dots,n)$
 is conjugate to a (non-trivial) reflection in $G$, which implies that
 the residue matrix of $\hat{B}^{(n)}\,{\rm d}x_n$ at $x_n=z_i^0(x')$ $(i=1,\dots,n)$
 is of rank one and diagonalizable.
 So is the residue matrix of $B^{(n)}\, {\rm d}x_n$.
 Consequently, we find that
 there exists an $(n\times n)$-matrix $T_0(x')$ such that
 \begin{equation}\label{equ:GqT_0}
 B^{(n)}=-(x_nI_n-T_0(x'))^{-1}B_{\infty}
 \end{equation}
 (at least) at any generic point~$x'$. Put
 \[
 h(x)=x_n^n-s_1(x')x_n^{n-1}+\cdots +(-1)^ns_n(x'),\qquad
 h(x)B^{(n)}=\sum_{i=0}^{n-1}D_i(x')x_n^i,
 \]
 where $D_{n-1}(x')=-B_{\infty}$.
 Then, by induction, we find
 \begin{equation*}
 D_{n-i}(x')=-\big(T_0(x')^{i-1}-s_1(x')T_0(x')^{i-2}+
 \dots +(-1)^{i-1} s_{i-1}(x')\big) B_{\infty},\qquad 1\le i\le n.
 \end{equation*}
 In particular, it holds that
 $T_0(x')=s_1(x')I_n-D_{n-2}(x')B_{\infty}^{-1}$,
 which implies that all entries of~$T_0(x')$ are weighted homogeneous
 polynomials in~$x'$.

 Finally we prove that the differential system (\ref{equ:w(T)=w(Mh)}) is normalized as~(\ref{eq:okubopfaff0}), that is,
 $B^{(k)}$ is written as
 \begin{equation} \label{eq:normal}
 B^{(k)}(x)=-(x_nI_n-T_0(x'))^{-1}\tilde{B}^{(k)}(x')B_{\infty}
 \end{equation}
 for some $\tilde{B}^{(k)}(x')\in \mathbb{C}[x']^{n\times n}$.
 We can write $h(x)B^{(k)}$ as a polynomial in $x_n$ to
 \begin{equation*}
 h(x)B^{(k)}=\sum_{j=0}^{m_k}D^{(k)}_j\,x_n^j,\qquad 1\le k\le n,
 \end{equation*}
 where $D^{(k)}_j \in \mathbb{C}[x']^{n\times n}$, $ D^{(n)}_{n-1}=-B_{\infty}$.
 From the equalities
 \[
 w\big(\big(B^{(k)}\big)_{i,j}\big)=w\big(\big(\hat{B}^{(k)}\big)_{i,j}\big)
 =w(x_j)-w(x_i)-w(x_k),
 \]
 we have $m_k\le n$ for $k\le n-1$ and $m_n=n-1$.
 If $m_k=n$, from the equality
 \mbox{$\big[D^{(k)}_n, -B_{\infty}\big]=O$},
 we have $\big(D^{(k)}_n\big)_{i,j}=0$ if $d_i\ne d_j$.
 If $d_i=d_j$, then $w\bigl(\big(D^{(k)}_n\big)_{i,j}\bigr) =-w(x_k)<0$,
 which concludes that $D^{(k)}_n=O$.
 Hence we have $m_k=n-1$,
 which implies~(\ref{eq:normal}) for $\tilde{B}^{(k)}(x')\in \mathbb{C}[x']^{n\times n}$.
 The Pfaffian system~(\ref{equ:w(T)=w(Mh)}) is evidently completely integrable
 because it admits a system of fundamental solutions.
 Then, in virtue of Lemma~\ref{lem:lemma2.2},
 $T(x)=T_0(x')-x_nI_n$, $\tilde{\Omega}=\sum_{k=1}^n\tilde{B}^{(k)}\,{\rm d}x_k$, $B_{\infty}$ satisfy the equations~(\ref{eq:commute}) and (\ref{eq:ci1}),
 from which we can conclude that~(\ref{equ:w(T)=w(Mh)}) is a reduced form of an extended Okubo system.
\end{proof}

\begin{Remark} \label{rem:gs vs okubo}
 The extended Okubo system (\ref{equ:w(T)=w(Mh)}) is called the $G$-{\it quotient system} (cf.~\cite{HK,KS}).
 The system of fundamental solutions to (\ref{equ:w(T)=w(Mh)}) is given by
 \[
 Y=-B_{\infty}^{-1}R(x')\big(\begin{matrix} V_n^{(h)} &V_{n-1}^{(h)} & \dots & V_1^{(h)} \end{matrix}\big)^t y(u)
 \]
 for any homogeneous linear function $y(u)$ of $u$.
 Therefore it is clear that the monodromy representation of the $G$-quotient system (\ref{equ:w(T)=w(Mh)}) coincides with
 the standard representation on $U_n$ of the well-generated complex reflection group $G$.
\end{Remark}

\begin{Theorem}\label{prop:C_1}
 There is a distinguished system of generators $\big\{F_i^{{\rm fl}}(u)\big\}_{1\leq i\leq n}$ of $G$-invariant polynomials $\mathbb{C}[u]^G$
 such that $F^{{\rm fl}}_{i}(u)$ is a homogeneous polynomial of degree $d_i$ and satisfies the following condition:
 let $t_i:=F^{{\rm fl}}_{i}(u)$ be regarded as a coordinate system on $X=U_n/G$.
 Then $t=(t_1,\dots,t_n)$ is a flat coordinate system on~$X$ induced by the $G$-quotient system~\eqref{equ:w(T)=w(Mh)}
 using Theorem~{\rm \ref{saitojacobian}}.

 Particularly, if $d_1<d_2<\dots <d_n$, then $F^{{\rm fl}}_{i}(u)$ are unique up to constant multiplications.
\end{Theorem}

\begin{proof}Put
\begin{equation} \label{eq:TBT0B}
 T(x)=T_0(x')-x_n I_n,
 \end{equation}
 where $T_0(x')$ is defined by (\ref{equ:GqT_0}).
 Let $C(x)$ be a matrix such that $T=-V_1\,C$.
Then (\ref{equ:w(T)=w(Mh)}) is written as
\begin{equation*}
 {\rm d}Y=T^{-1}\,{\rm d}C\,B_{\infty} Y
 =-(V_1 C)^{-1}\,{\rm d}C\,B_{\infty} Y.
\end{equation*}

Put
\begin{equation*}
 F^{{\rm fl}}_j(u):=C_{nj}(F_1(u),F_2(u),\dots,F_n(u)),\qquad 1\le j\le n.
\end{equation*}
Since $w(C_{nj})=w(x_j)$
and $R(x')$ has the form~(\ref{eq:R}),
the equalities
$t_j=C_{nj}(x)$, $1\le j\le n$
are solved
by weighted homogeneous polynomials $x_j=x_j(t)$, $1\le j\le n$,
and hence $\big\{ F^{{\rm fl}}_j(u)\big\}$ is a fundamental system of $G$-invariant polynomials.
In virtue of Theorem~\ref{saitojacobian},
$(t_1,\dots,t_n)=\big(F_1^{{\rm fl}}(u),\dots,F_n^{{\rm fl}}(u)\big)$ is the flat coordinate system on~$X$ induced
by the $G$-quotient system~(\ref{equ:w(T)=w(Mh)}).

The uniqueness of $F_i^{{\rm fl}}(u)$, $1\leq i \leq n$ up to constant multiplications under the condition
 $d_1 < d_2 < \cdots < d_n$
is clear from Remarks~\ref{rem:flat pvf constant} and~\ref{rem:uniqunessofSaitoOkubo}.
\end{proof}

\begin{Definition} \label{def:flatbasicinv}
 The system of generators $\big\{F^{{\rm fl}}_{i}(u)\big\}$ in Theorem~\ref{prop:C_1}
 is called a {\it flat generator system} of $G$-invariant polynomials or {\it flat basic invariants} of~$G$.
\end{Definition}

\begin{Remark} \label{rem:G-quotflat}
 Take the flat coordinate system $t=(t_1,\dots,t_n)$ on $X=U_n/G$ in Theorem~\ref{prop:C_1}.
 Let~$M_{\mathcal{V}^{(h)}}(t)$ and $V_i^{(h)}$ be defined by (\ref{equ:Def of M_V and M_h}) and (\ref{equ:LogVFs}) respectively
 with respect to $h=h(t)$.
 Let~$C(t)$ be the matrix whose entries are weighted homogeneous polynomials satisfying $V_1^{(h)}C(t)\allowbreak =M_{\mathcal{V}^{(h)}}(t)$.
 Then the $G$-quotient system (\ref{equ:w(T)=w(Mh)}) is written as
 \begin{equation} \label{eq:G-quotflat}
 {\rm d}Y=\big[{-}\big(V_1^{(h)}C(t)\big)^{-1}{\rm d}C(t)\,B_{\infty}\big]Y
 \end{equation}
 in the flat coordinate~$t$.
 The system of fundamental solutions to~(\ref{eq:G-quotflat}) is given by
 \begin{equation*} 
 Y=-B_{\infty}^{-1}\big(\begin{matrix} V_n^{(h)} & V_{n-1}^{(h)} & \dots & V_1^{(h)} \end{matrix}\big)^t y(u)
 \end{equation*}
 for any homogeneous linear function $y(u)$ of $u$.
\end{Remark}

\begin{Corollary} \label{cor:polyG}
 The potential vector field $\vec{g}=(g_1,\dots,g_n)$ for the $G$-quotient system \eqref{eq:G-quotflat}
 has polynomial entries in $t=(t_1,\dots,t_n)$.
\end{Corollary}

\begin{proof} It is clear that $h(t)$ and the entries of $M_{\mathcal{V}^{(h)}}(t)$ are polynomials in~$t$.
 Thus, the entries of $C(t)$ in Remark~\ref{rem:G-quotflat} are also polynomials in~$t$.
 Then by~(\ref{eq:pvfC}) in Remark~\ref{rem:pvfC}, we see that the entries of~$\vec{g}$ are polynomials in~$t$.
\end{proof}

\begin{Example} \label{ex:G26} Taking $G=G_{26}$ which is the complex reflection group of No.~26 in the table of~\cite{ST},
 we explain how to construct the flat coordinate system and the potential vector field by our method.
We start with the classical basic $G$-invariants given in~\cite{Mas}:
\begin{gather*}
x_1=F_1(u)=
u_1^6 - 10 u_1^3 u_2^3 - 10 u_1^3 u_3^3 + u_2^6 - 10 u_2^3 u_3^3 + u_3^6,\\
x_2=F_2(u)=\big(u_1^3+u_2^3+u_3^3\big) \big(\big(u_1^3+u_2^3+u_3^3\big)^3+216 (u_1 u_2 u_3)^3\big),\\
x_3=F_3(u)=\big(\big(u_2^3-u_3^3\big) \big(u_3^3-u_1^3\big) \big(u_1^3-u_2^3\big)\big)^2.
\end{gather*}
The discriminant $h(x)$ and the corresponding Saito matrix $M_{\mathcal{V}}=M_{\mathcal{V}^{(h)} }$
defined in~(\ref{equ:Def of M_V and M_h}) are given by
\begin{gather}
h=\big(1/432^2\big) x_3 \big(\big(x_1^3-3 x_1 x_2-432 x_3\big)^2-4 x_2^3\big), \label{eq:h_x} \\
M_{\mathcal{V}^{(h)}}=
\begin{pmatrix}
 \dfrac{ - x_1^3 + x_1 x_2 + 432 x_3}{432}&
 \dfrac{x_2 ( - x_1^2 + x_2)}{216}& 0\vspace{1mm}\\
 \dfrac{x_2}{216}& \dfrac{ - x_1^3 + 5 x_1 x_2 + 432 x_3}{432}&0\vspace{1mm}\\
 \dfrac{x_1}{3}&\dfrac{2 x_2}{3}&x_3
 \end{pmatrix}.\nonumber
\end{gather}

The matrices $P_k:=\big(\tilde{V}_k M_{\tilde{\mathcal{V}}}\big) M_{\tilde{\mathcal{V}}}^{-1}$
in Lemma~\ref{lemma:tildeVi is log} (ii) are computed as
\begin{gather*}
P_3 =
\begin{pmatrix}
\dfrac{5 ( - x_1^2 + x_2)}{2592}&\dfrac{ -x_3}{3}&
\dfrac{x_1^4 - 2 x_1^2 x_2 - 432 x_1 x_3 + x_2^2}{46656}\vspace{1mm}\\
 \dfrac{ -x_1}{648}&\dfrac{7 ( - x_1^2 + x_2)}{2592}&
 \dfrac{ - x_1^3 + x_1 x_2 + 432 x_3}{93312}\vspace{1mm}\\
 \dfrac{1}{18}&0&0
\end{pmatrix},\\
P_2 =
\begin{pmatrix}
 \dfrac{x_1}{1296} &\dfrac{ - x_1^2 + x_2}{2592}&
 \dfrac{ - x_1^3 + x_1 x_2 + 432 x_3}{93312}\vspace{1mm}\\
 \dfrac{7}{2592}&\dfrac{5 x_1}{1296}&\dfrac{x_2}{46656}\vspace{1mm}\\
 0& {1\over 18}&0
\end{pmatrix},\\
P_1 =\operatorname{diag}\left[{13\over 18},{7\over 18},{1\over 18}\right]
=-B_{\infty},
\end{gather*}
and the coefficients $\hat{B}^{(k)}$ of the Pfaffian system~(\ref{equ:G-quotient}) are obtained by
$\hat{B}^{(k)}=\sum_{j=1}^3\big(M_{\mathcal{V}}^{-1}\big)_{k,j} P_{4-j}$.

The matrix $\hat{B}_{\infty}(x')=-(x_3\hat{B}^{(3)})|_{x_3=\infty}$ in the proof of Lemma~\ref{lemma:w(T)=w(Mh)} is
\[
\hat{B}_{\infty}(x')=
\begin{pmatrix}
 \dfrac{13}{18}&\dfrac{x_1}{9}&\dfrac{x_1^2-x_2}{324}\vspace{1mm}\\
 0&\dfrac{7}{18}&\dfrac{-x_1}{648}\vspace{1mm}\\
 0&0&\dfrac{1}{18}
 \end{pmatrix},
\]
 which is diagonalized by
\[ R(x')=
\begin{pmatrix}
 1 &\dfrac{x_1}{3}&\dfrac{5x_1^2-6 x_2}{1296}\vspace{1mm}\\
 0&1&\dfrac{-x_1}{216}\vspace{1mm}\\
 0&0&1
 \end{pmatrix}
 \]
 as
 \[ R(x') \hat{B}_{\infty}(x') R(x')^{-1}
 =\operatorname{diag}\left[{1\over 3},{2\over 3},1\right]-{19\over 18}I_3
 =B_{\infty}.
 \]
The coefficient matrix $B^{(3)}(x)$ in the $x_3$-direction of (\ref{equ:w(T)=w(Mh)})
in Lemma~\ref{lemma:w(T)=w(Mh)} is computed as
\[
 B^{(3)}(x) = B_{\infty}^{-1} R(x') \hat{B}^{(3)} R(x')^{-1}B_{\infty}.
\]
Using (\ref{equ:GqT_0}), we obtain $T(x)$ in (\ref{eq:TBT0B}) as
\[
 T(x)=B_{\infty}\big(B^{(3)}(x)\big)^{-1}.
\]
Let $C_{ij}=-(1-w(x_i)+w(x_j))^{-1}T_{ij}$, $1\leq i,j\leq 3$.
Then the matrix $C=(C_{ij})_{1\leq i,j\leq 3}$ is given by
\begin{equation*}
 C=
 \begin{pmatrix}
 \dfrac{-4x_1^3 + 9x_1x_2+3888x_3}{3888} & \dfrac{-5x_1^4+12x_1^2x_2+18x_2^2}{15552} & \dfrac{x_1(2x_1^4-9x_1^2x_2+18x_2^2)}{1679616} \vspace{1mm}\\
 \dfrac{-x_1^2+3x_2}{432} & \dfrac{-7x_1^3+27x_1x_2+3888x_3}{3888} & \dfrac{-x_1^4+3x_1^2x_2+9x_2^2}{559872} \\
 x_1 & \dfrac{-x_1^2+6x_2}{6} & \dfrac{-7x_1^3+18x_1x_2+3888x_3}{3888}
 \end{pmatrix}.
\end{equation*}
For any constant number $c_i\neq 0$, $t_i=c_i C_{3,i}$, $1\leq i\leq 3$ give the flat basic invariants,
and
\begin{equation*}
 g_i=\frac{c_ic_3}{1+w(x_i)}\sum_{j=1}^3w(x_j)C_{j,i}C_{3,j},\qquad 1\leq i\leq 3
\end{equation*}
give the potential vector field $\vec{g}_{26}=(g_1,g_2,g_3)$ (see Remarks~\ref{rem:pvfC} and~\ref{rem:flat pvf constant}).
If we choose $c_1=1/6^{4/3}$, $c_2=-1/\big(2\cdot 6^{5/3}\big)$, $c_3=1$,
then we obtain the flat basic invariants
\begin{gather}
 t_1= 6^{-4/3} x_1, \qquad
 t_2={1\over 2\cdot 6^{8/3}}\big(x_1^2-6x_2\big), \label{eq:t1t2} \\
 t_3= x_3 +{1\over 6^3} x_1 x_2 - {7\over 3\cdot 6^4} x_1^3, \label{eq:t3}
\end{gather}
and the potential vector field $\vec{g}_{26}=(g_1,g_2,g_3)$ with
\begin{gather}
 g_1= t_1 t_3+{1\over 2}t_2^2+{1\over 2}t_1^2 t_2+{1\over 8}t_1^4,\qquad
 g_2= t_2 t_3-{1\over 2}t_1 t_2^2+{1\over 2}t_1^3 t_2+{1\over 8}t_1^5, \label{eq:g1g2} \\
 g_3= {1\over 2}t_3^2-{1\over 3}t_2^3+t_1^2 t_2^2+{1\over 4}t_1^4 t_2+{1\over 6}t_1^6. \label{eq:g3}
\end{gather}
\end{Example}

\begin{Remark} \label{rem:KMS}
 Irreducible finite complex reflection groups are classified by Shephard--Todd~\cite{ST}:
 except for rank 1 groups, there are two infinite families $A_n$, $G(pq,p,n)$,
 plus~34 exceptional groups $G_4,G_5,\dots,G_{37}$.
 Explicit formulas of the potential vector fields for the exceptional groups $G_{24},G_{27},G_{29},G_{33},G_{34}$
 are found in~\cite{KMS}.
 The results in~\cite{KMS} are obtained by solving extended WDVV equations related with the
discriminants of the groups given above.
It is underlined that the potential vector fields
given there are shown to be
unique up to constant multiplications of the variables
and therefore they correspond to $G$-quotient systems for these groups.\footnote{The sentence in the last two lines of \cite[p.~19]{KMS} is not correct.
Namely, ``if $G$ is an irreducible finite complex reflection groups
and it is well-generated in the sense of~\cite{AL0},
there is a unique flat structure for the discriminant of $G$'' is not correct,
cf.\ Remark~\ref{rem:ShephardFrobenius}.}
 As explained there in the case $G_{34}$, the potential vector field given there
was not proved to be that for $G_{34}$ (i.e., it was a conjecture).
 Later the authors confirmed the result affirmatively.
This is shown with the aid of the result by Terao--Enta \cite[Appendix]{Or}
 and the computation on the Saito matrix of~$G_{34}$ case by D.~Bessis and J.~Michel (private communication with J.~Michel).

 It is possible to compute the explicit form of flat basic invariants for
 a complex reflection group~$G$ from its potential vector field by the following procedure:
 let $(x_1,\dots,x_n)=(F_1(u),\dots,F_n(u))$ be arbitrary basic invariants of $G$
 and write down the discriminant $h_x(x)$ of $G$ in terms of $x$.
 (In order to avoid confusion, we denote by $h_x(x)$ the discriminant $h$ in the coordinate $x$.)
 On the other hand, one can write down the discriminant $h_t(t)$ of $G$ in terms of $t$ from the potential vector field
 in the manner described in Section~\ref{sec:saitostrandOkubo}
 (i.e., $h_t(t)=\det (E\mathcal{C})=\det \bigl((1+w_j-w_i)\frac{\partial g_j}{\partial t_i}\bigr)$).
 Find a weight preserving coordinate change $t=t(x)$ such that $h_t(t(x))=h_x(x)$.
 Then $\bigl(F_1^{{\rm fl}}(u),\dots,F_n^{{\rm fl}}(u)\bigr):=
 \bigl(t_1(F_1(u),\dots,F_n(u)),\dots,t_n(F_1(u),\dots,F_n(u))\bigr)$ gives flat basic invariants of $G$.
 Let us demonstrate this procedure for the example $G_{26}$.
 By~(\ref{eq:h_x}), we have
 \begin{equation*}
 h_x(x)=\big(1/432^2\big) x_3 \big(\big(x_1^3-3 x_1 x_2-432 x_3\big)^2-4 x_2^3\big).
 \end{equation*}
 On the other hand, from the potential vector field $\vec{g}_{26}$ of $G_{26}$
 given by~(\ref{eq:g1g2}),~(\ref{eq:g3}), we have
 \begin{gather*}
 h_t(t)= -\frac{1}{108}\big(3t_3+4t_1^3+6t_1t_2\big)\big(5t_1^6-48t_1^4t_2+12t_3t_1^3+36t_1^2t_2^2+72t_3t_1t_2-36t_3^2-48t_2^3\big).
 \end{gather*}
 Any weight preserving coordinate change $t=t(x)$ (satisfying the condition~(\ref{cond:T}))
 is written as
 \begin{equation*}
 t_1=c_1x_1,\qquad t_2=c_2x_2+c_3x_1^2, \qquad t_3=x_3+c_4x_2x_1+c_5x_1^3,
 \end{equation*}
 where we regard $c_1,\dots,c_5$ as unknown coefficients.
 Solve $h_t(t(x))=h_x(x)$ for $c_1,\dots,c_5$.
 Then we have
 \begin{gather*}
 c_1=6^{-4/3},\qquad\!\! c_2=-1/\big(2\cdot 6^{5/3}\big),\qquad\!\! c_3=1/\big(2\cdot 6^{8/3}\big),\qquad\!\! c_4=6^{-3},\qquad\!\! c_5=-7/\big(3\cdot 6^4\big),
 \end{gather*}
 which coincide with (\ref{eq:t1t2}),~(\ref{eq:t3}).
\end{Remark}

\begin{Remark}
 After a preprint of the present paper had been released on arXiv,
 Arsie--Lorenzoni \cite{AL2} explicitly computed (by case-by-case analysis)
 examples of the flat basic invariants for rank 2 and 3 well-generated complex reflection groups
 based on their theory of bi-flat $F$-manifolds
 (it was proved in \cite{AL2,KoMiSh} that Arsie--Lorenzoni's bi-flat $F$-manifold
 is equivalent to Sabbah's Saito structure).
 Recently Konishi--Minabe--Shiraishi \cite{KoMiSh} studied the Saito structure
 constructed in Theorem~\ref{prop:C_1}
 from the viewpoint of an extension of Dubrovin's almost duality \cite{Du2}.
 The Saito structure induced by the $G$-quotient system in Theorem~\ref{prop:C_1} is same as
 the ``standard generalized Saito flat coordinates'' in \cite{AL2}
 and the ``natural Saito structure for a complex reflection group'' in \cite{KoMiSh}.
 Indeed the ``dual connection'' adopted in \cite{AL2,KoMiSh} is nothing but the meromorphic connection whose horizontal sections
 satisfy an extended Okubo system.
 This relationship was clarified in \cite{KoMiSh}.
\end{Remark}

\begin{Remark} \label{rem:ShephardFrobenius}
 Shephard groups are the symmetry groups of regular complex polytopes,
 which form a~subclass of complex reflection groups.
 The papers \cite{Du2, OS} treat Frobenius structures constructed on the orbit spaces of Shephard groups.
 We shall compare the results in \cite{Du2, OS} with Theorem~\ref{prop:C_1} in the present paper.
 By Theorem~\ref{prop:C_1}, we see that the Saito structures on the orbit spaces induced by the $G$-quotient systems for Shephard groups
 are divided into the following two types:
 \begin{enumerate}\itemsep=0pt
 \item[(i)] the Saito structure has only a potential vector field (but has no prepotential),
 \item[(ii)] the Saito structure has a prepotential, that is the Saito structure becomes a Frobenius structure.
 \end{enumerate}
 (The existence of the type (i) is discussed also in~\cite{AL2}.)
 It is known that to each Shephard group $G$ there is an associated Coxeter group $W$
 whose discriminant is isomorphic to that of the Shephard group~$G$.
 Then we see that there are (at least) two extended Okubo systems on the orbit space of the Shephard group $G$
 which have singularities along the discriminant:
 one is the $G$-quotient system of the Shephard group~$G$,
 the other is the $G$-quotient system of the associated Coxeter group~$W$.
 The Frobenius structure on the orbit space of a Shephard group~$G$ described in \cite{Du2, OS} corresponds to the $G$-quotient system
 of the associated Coxeter group~$W$.

 $G_{26}$ in Example~\ref{ex:G26} is a Shephard group belonging to the type~(i).
As is pointed out in~\cite{AL2}, the potential vector field $\vec{g}_{26}$ of $G_{26}$ given by (\ref{eq:g1g2}),~(\ref{eq:g3})
admits a one-parameter deformation~$\vec{g}(q)$ (where $q\in\mathbb{C}\setminus\{0\}$) given by
\begin{gather*}
g_1= t_1 t_3 + {(q - 2) (q - 3)\over 6} t_1^4 + ( - q + 2) t_1^2 t_2
+{1\over 2} t_2^2,\\
g_2=
 t_2 t_3-{(q^2 - 1) (q - 2)\over 5}\,t_1^5
+ {2 q (q - 1)\over 3}t_1^3 t_2 +( - q + 1) t_1 t_2^2,\\
g_3=
 {1\over 2}t_3^2
 +{4 (q^2 - 2 q + 2) (q - 1)\over 15}t_1^6 -(q - 1) (q - 2) t_1^4 t_2 + 2(q - 1) t_1^2 t_2^2 -{1\over 3} t_2^3,
\end{gather*}
such that $\vec{g}(3/2)=\vec{g}_{26}$.
For $q=2$, we notice that $\vec{g}(2)$ has a prepotential
\[
 F(t)=\frac{t_1 t_3^2 + t_2^2 t_3}{2}+{8\over 105} t_1^7 + {2\over 3} t_1^3 t_2^2 -{1\over 3} t_1 t_2^3
\]
satisfying
\[
 g_1={\partial F\over \partial t_3},\qquad
g_2={\partial F\over \partial t_2},\qquad
g_3={\partial F\over \partial t_1},
\]
which coincides with the prepotential for the Weyl group of type~$B_3$.
The Weyl group $B_3$ is nothing but the Coxeter group associated with the Shephard group~$G_{26}$.

In general, an extended Okubo system is constructed from a potential vector field $\vec{g}$
and a~diagonal matrix $B_{\infty}=\operatorname{diag}[w_1-\lambda,w_2-\lambda,\dots, w_n-\lambda]$
(where $\lambda$ is the same parameter as in Remark~\ref{rm:henkan}):
\begin{equation} \label{eq:G-quotflat2}
 {\rm d}Y=\big[{-}\big(EC(t)\big)^{-1}\,{\rm d}C(t)\,B_{\infty}\big]Y, \qquad \mbox{for} \quad C(t):=\left(\frac{\partial g_j}{\partial t_i}\right).
\end{equation}
In this example, $n=3$, $\vec{g}=\vec{g}(q)$ and $w_1=1/3, w_2=2/3, w_3=1$.
Then the local exponents of the Okubo system~(\ref{eq:G-quotflat2}) at $t_3=z_i^0$, $i=1,2,3$ are given by
\[
 (0,0,1/q+c),\qquad (0,0,1/2+c),\qquad (0,0,1/2+c),
\]
where $z_i^0$, $1\leq i\leq 3$ are defined by $h(t)=\big(t_3-z_1^0\big)\big(t_3-z_2^0\big)\big(t_3-z_3^0\big)$ and $c:=-1+\lambda -{1\over 3 q}$.
This implies that the monodromy group of the extended Okubo system~(\ref{eq:G-quotflat2}) varies depending on the parameters $q$ and $\lambda$.
For $q=3/2$ and $\lambda=19/18(=1+1/d_3)$ (cf.~(\ref{equ:B_inf})),
 the corresponding extended Okubo system is nothing but
the $G$-quotient system for $G_{26}$, and its monodromy group is isomonodromic to $G_{26}$.
For $q=2$ and $\lambda=7/6$,
the extended Okubo system~(\ref{eq:G-quotflat2}) is the $G$-quotient system for $B_3$,
and its monodromy group is isomorphic to the Weyl group of type~$B_3$.
In this case, the Saito structure in Theorem~\ref{prop:C_1} for~$G_{26}$ is distinct to the Frobenius structure constructed in~\cite{Du2, OS}
on the orbit space of~$G_{26}$.
We remark that the potential vector field~$\vec{g}(q)$ corresponds to a one-parameter algebraic solution to the Painlev\'e VI equation
 (Solution~III in~\cite{LT}).
See~\cite{KMS4} for details on Solution III and its potential vector field.

We take a Shephard group $G_{32}$ as another example.
The Saito structure induced by the $G$-quotient system of $G_{32}$ in Theorem~\ref{prop:C_1} has a prepotential given by
\begin{equation*}
 F_{32}:=\frac{2t_1^6}{15}+\frac{2t_1^3t_2^2}{3}+t_1^2t_3^2+t_1t_2^2t_3+\frac{t_1t_4^2}{2}
 +\frac{t_2^4}{12}+t_2t_3t_4+\frac{t_3^3}{3},
\end{equation*}
which implies that $G_{32}$ belongs to the type (ii).
$F_{32}$ coincides with the prepotential of the Frobenius structure for the Weyl group of type~$A_4$
up to constant multiplications of the flat coordinates.
The Weyl group of type $A_4$ is the associated Coxeter group with the Shephard group~$G_{32}$.
Put $\vec{g}_{32}:=(\partial F_{32}/\partial t_4,\partial F_{32}/\partial t_3,\partial F_{32}/\partial t_2,\partial F_{32}/\partial t_1)$.
Then, for $\vec{g}=\vec{g}_{32}$ and \mbox{$\lambda=31/30$},
the extended Okubo system~(\ref{eq:G-quotflat2}) is the $G$-quotient system of $G_{32}$ and its monodromy group is isomorphic to~$G_{32}$.
For $\vec{g}=\vec{g}_{32}$ and $\lambda=6/5$, (\ref{eq:G-quotflat2}) is the $G$-quotient system for $A_4$ and its monodromy group is
isomorphic to the Weyl group of type $A_4$.
In this case, the Saito structure in Theorem~\ref{prop:C_1} for $G_{32}$ is identical with the Frobenius structure
constructed on the orbit space of $G_{32}$ in~\cite{Du2, OS}.
\end{Remark}

\section[Examples of potential vector fields corresponding to algebraic solutions to the Painlev\'e VI equation]{Examples of potential vector fields corresponding\\ to algebraic solutions to the Painlev\'e VI equation} \label{sec:algsol}

In this section we show some examples of potential vector fields in three variables
which correspond to algebraic solutions to the Painlev\'e~VI equation.

\looseness=1 Algebraic solutions to the Painlev\'e~VI equation were studied and constructed by
 many authors including N.J.~Hitchin \cite{Hit1,Hit2},
 B.~Dubrovin~\cite{Du}, B.~Dubrovin and M.~Mazzocco~\cite{DM},
P.~Boalch \cite{BoC,Bo0,Bo1,Bo3,Bo2},
A.V.~Kitaev~\cite{Kit1,Kit2}, A.V.~Kitaev and R.~Vid\=unas~\cite{KV,VK}, K.~Iwasaki~\cite{Iwa}.
In his review article~\cite{Bo4}, Boalch classified all the known algebraic solutions to the Painlev\'e VI
based on types of the monodromy groups of the associated $2\times 2$ or $3\times 3$ linear differential equations.
After that, by classifying finite orbits of the extended modular group action on conjugacy classes of ${\rm SL}(2,\mathbb{C})$-triples,
Lisovyy and Tykhyy~\cite{LT} proved that Boalch's list was complete, that is,
they confirmed that there is no other solution than those in Boalch's list.
One of the principal aims of our study is the determination of a flat coordinate system
and a potential vector field for each of such algebraic solutions.
In spite that this aim is still not succeeded because of
complexity of computation, we show some examples of potential vector fields.
Some of the results below are already given in~\cite{KMS}.
Other examples can be found in~\cite{KMS4}.
From the construction of polynomial potential vector fields corresponding to finite complex reflection groups of rank three
(Corollary~\ref{cor:polyG}) and Corollary~\ref{cor:corgexWDVVandPVI},
we obtain a class of algebraic solutions to the Painlev\'e~VI equation.
The relationship between finite complex reflection group of rank three
and solutions to the Painlev\'e VI equation was first studied by Boalch~\cite{BoC}.
(More precisely speaking, it was conjectured in \cite{BoC} that the solutions obtained
from finite complex reflection groups by his construction are algebraic
and this conjecture was proved in his succeeding papers.)
The construction in the present paper answers the question~3) in the last part of~\cite{BoC}:
``Is there a geometrical or physical interpretation of these solutions?''

It is known by Hertling \cite{He} that any polynomial prepotential is obtained from
a Frobenius manifold on the orbit space of a (product of) finite Coxeter group.
Surprisingly enough, there are examples in Sections~\ref{section6.4}--\ref{section6.6}, whose potential vector fields have polynomial entries but the corresponding Saito structures are not isomorphic to one on the orbit space of any finite complex reflection group
because the free divisors defined by $F_{B_{6}}$, $F_{H_2}$, $F_{E_{14}}$ in Sections~\ref{section6.4}--\ref{section6.6} respectively are not isomorphic to the discriminant of any finite complex reflection group.
The existence of these examples suggests that an analogue of Hertling's theorem does not hold in the case of non Frobenius manifolds.
It seems to be an interesting problem to classify polynomial potential vector fields.

To avoid confusion, we prepare the convention which will be used in the following.
We treat the case $n=3$.
Let $t=(t_1,t_2,t_3)$ be a flat coordinate system and $\vec{g}=(g_1,g_2,g_3)$ denotes a~potential vector field.
Let $w(t_i)$ be the weight of $t_i$ and assume $0<w(t_1)<w(t_2)<w(t_3)=1$.
The matrix $C$ is defined by $C=\big(\frac{\partial g_j}{\partial t_i}\big)$
and $T=-EC=-\sum_{j=1}^3w(t_j)t_j\partial_{t_j}C$.
In this section, an algebraic solution LT$n$ stands for ``Solution~$n$'' in Lisovyy--Tykhyy \cite[pp.~156--162]{LT},
and~$G_n$ denotes the finite complex reflection group of No.~$n$ in the table of Shephard--Todd~\cite{ST}.

\subsection{Algebraic solutions related with icosahedron}

We treat the three algebraic solutions to Painlev\'e VI obtained by Dubrovin~\cite{Du} and Dubrovin--Mazzocco~\cite{DM}.
The solutions in this subsection have prepotentials.

{\bf Icosahedral solution $\boldsymbol{(H_3)}$ (icosahedral solution 31 \cite{Bo1,Bo4}, LT16).} In this case,
\[
w(t_1)=\frac{1}{5},\qquad w(t_2)=\frac{3}{5},\qquad w(t_3)=1
\]
and there is a prepotential defined by
\begin{equation*}
F=\frac{ t_2^2t_3 +t_1t_3^2}{2}+\frac{t_1^{11}}{3960} +\frac{ t_1^5t_2^2}{20} + \frac{t_1^2t_2^3}{6}.
\end{equation*}

We don't enter the details on this case. See \cite{Du,DM}.

{\bf Great icosahedral solution $\boldsymbol{(H_3)'}$ (icosahedral solution~32 \cite{Bo1,Bo4}, LT17).}
Let $(t_1,t_2,t_3)$ be a flat coordinate system and their weights are given by
\[
w(t_1)=\frac{3}{5},\qquad w(t_2)=\frac{4}{5},\qquad w(t_3)=1.
\]

We introduce an algebraic function $z$ of $t_1$, $t_2$ defined by the relation
\[
t_2 + t_1z + z^4=0.
\]
It is clear from the definition that $w(z)=\frac{1}{5}$.
In this case, the prepotential is an algebraic function of $(t_1,t_2,t_3)$
defined by
\[
F=\frac{t_2^2t_3 + t_1t_3^2}{2} - \frac{t_1^4z}{18} -
\frac{7t_1^3z^4}{72} - \frac{17t_1^2z^7}{105} - \frac{2t_1z^{10}}{9} - \frac{64z^{13}}{585}.
\]

{\bf Great dodecahedron solution $\boldsymbol{(H_3)''}$ (icosahedral solution~41 \cite{Bo1,Bo4}, LT31).}
Let $(t_1,t_2,t_3)$ be a flat coordinate system and their weights are given by
\[
w(t_1)=\frac{1}{3},\qquad w(t_2)=\frac{2}{3},\qquad w(t_3)=1.
\]

We introduce an algebraic function $z$ of $t_1$, $t_2$ defined by the relation
\[
 -t_1^2 + t_2 + z^2=0.
\]
It is clear from the definition that $w(z)=\frac{1}{3}$.
In this case, the prepotential is an algebraic function of $(t_1,t_2,t_3)$
defined by
\[
F= \frac{t_2^2t_3 + t_1t_3^2}{2}+\frac{4063t_1^7}{1701} +
\frac{19t_1^5z^2}{135} - \frac{73t_1^3z^4}{27} + \frac{11t_1z^6}{9} - \frac{16z^7}{35}.
\]

\begin{Remark}
The prepotential $F$ for the icosahedral solution $(H_3)$ was firstly obtained by B.~Dubrovin~\cite{Du}.
The above algebraic solutions including the remaining two cases $(H_3)'$, $(H_3)''$ were treated by B.~Dubrovin and M.~Mazzocco~\cite{DM}.
The authors were informed by B.~Dubrovin that his student Alejo Keuroghlanian computed the algebraic Frobenius manifold
for the case of the great icosahedron $(H_3)'$ in his master thesis
``Varieta di Frobenius algebraiche di dimensione~$3$''.
The authors don't know whether
the potential for $(H_3)''$ is known or not.
Topics on these solutions are treated in~\cite{KMS3}.
\end{Remark}

\subsection[Algebraic solution related with the complex reflection group $G_{24}$ \cite{Bo0}]{Algebraic solution related with the complex reflection group $\boldsymbol{G_{24}}$ \cite{Bo0}}

The following case is related with the complex reflection group~$G_{24}$.

{\bf Klein solution of Boalch \cite{Bo0,Bo4} (LT8).}
In this case,
\begin{gather*}
w(t_1)=\frac{2}{7},\qquad w(t_2)=\frac{3}{7},\qquad w(t_3)=1,\\
g_1 = \big({-}2t_1^3t_2 + t_2^3 + 12t_1t_3\big)/12,\\
g_2 = \big(2t_1^5 + 5t_1^2t_2^2 + 10t_2t_3\big)/10, \\
g_3 = \big({-}8t_1^7 + 21t_1^4t_2^2 + 7t_1t_2^4 + 28t_3^2\big)/56.
\end{gather*}
The determinant $\det(-T)$ is regarded as the discriminant of the complex reflection group~$G_{24}$ if
$t_1$, $t_2$, $t_3$ are taken as basic invariants (cf.~\cite{ST}).

\subsection[Algebraic solutions related with the complex reflection group $G_{27}$]{Algebraic solutions related with the complex reflection group $\boldsymbol{G_{27}}$}

The following two cases are related with the complex reflection group $G_{27}$ (cf.~\cite{ST}).

{\bf Icosahedral solution 38 of Boalch \cite{Bo1,Bo4} (LT26).}
In this case,
\begin{gather*}
w(t_1)=\frac{1}{5},\qquad w(t_2)=\frac{2}{5},\qquad w(t_3)=1,\\
g_1 = \big({-}t_1^6 - 15t_1^4t_2 + 15t_1^2t_2^2 + 10t_2^3 + 30t_1t_3\big)/30,\\
g_2 = \big(5t_1^7 + 3t_1^5t_2 + 15t_1^3t_2^2 - 5t_1t_2^3 + 6t_2t_3\big)/6, \\
g_3 = \big({-}105t_1^{10} + 200t_1^8t_2 + 350t_1^6t_2^2 + 175t_1^2t_2^4 - 14t_2^5 + 20t_3^2\big)/40.
\end{gather*}
The determinant $\det(-T)$ is regarded as the discriminant of the complex reflection group~$G_{27}$,
if $t_1$, $t_2$, $t_3$ are taken as basic invariants.

{\bf Icosahedral solution 37 of Boalch \cite{Bo1,Bo4} (LT27).}
In this case, $z$ is an algebraic function of $t_1$, $t_2$ defined by
\begin{gather*}
 -t_2 - t_1z + 2z^3=0,\\
w(t_1)=\frac{2}{5},\qquad w(t_2)=\frac{3}{5},\qquad w(t_3)=1,\qquad w(z)=\frac{1}{5},
\\
g_1 = \big(175t_1t_3 - 70t_1^3z + 70t_1^2z^3 + 378t_1z^5 - 540z^7\big)/175, \\
g_2 = \big(10t_1^4 - 120t_1t_2^2 + 75t_2t_3 + 30t_1^2z^4 - 192t_1z^6 + 324z^8\big)/75,\\
g_3 = \big(16t_1^5 + 80t_1^2t_2^2 + 25t_3^2 - 80t_1^3z^4 + 540t_1^2z^6 - 1080t_1z^8 + 432z^{10}\big)/50.
\end{gather*}
The determinant $\det(-T)$ is regarded as the discriminant of the complex reflection group $G_{27}$,
if $z$, $t_1$, $t_3$ are taken as basic invariants.

\subsection[Algebraic solutions related with the polynomial $F_{B_6}$ in \cite{Se1}]{Algebraic solutions related with the polynomial $\boldsymbol{F_{B_6}}$ in \cite{Se1}}\label{section6.4}

We recall the polynomial
\[
F_{B_6}=9xy^4+6x^2y^2z-4y^3z+x^3z^2-12xyz^2+4z^3,
\]
which is a defining equation of a free divisor in $\mathbb{C}^3$ (cf.~\cite{Se1}).
There are two algebraic solutions which are related with the polynomial $F_{B_6}$.

{\bf Icosahedral solution 27 of Boalch \cite{Bo1,Bo4} (LT13).}
In this case,
\begin{gather*}
w(t_1)=\frac{1}{15},\qquad w(t_2)=\frac{1}{3},\qquad w(t_3)=1,
\\
g_1 = -\frac{1}{33}t_1\big(3t_1^{10}t_2 + 11t_2^3 - 33t_3\big),\\
g_2 = \frac{1}{76}\big({-}5t_1^{20} + 114t_1^{10}t_2^2 + 19t_2^4 + 76t_2t_3\big), \\
g_3 = \frac{1}{870}\big(100t_1^{30} + 1740t_1^{20}t_2^2 - 5220t_1^{10}t_2^4 + 116t_2^6 + 435t_3^2\big).
\end{gather*}
The determinant $\det(-T)$ regarded as a polynomial of $t_2$, $t_1^{10}$, $t_3$ coincides with $F_{B_6}$
by a weight preserving coordinate change up to a non-zero constant factor.

{\bf Solution obtained by Kitaev \cite{Kit2} (Icosahedral solution~26 \cite{Bo1,Bo4}, LT14).} In this case, $z$ is an algebraic function of $t_1$, $t_2$ defined by
\begin{gather*}
 t_1^2 + t_2z^6 + z^{16}=0,
\\
w(t_1)=\frac{8}{15},\qquad w(t_2)=\frac{2}{3},\qquad w(t_3)=1,\qquad w(z)=\frac{1}{15},\\
g_1 = -\big(2093t_1^4 - 897t_1t_3z^9 + 3450t_1^2z^{16} + 525z^{32}\big)/\big(897z^9\big), \\
g_2 = \big({-}238t_1^5 + 85t_2t_3z^{15} + 1700t_1^3z^{16} - 750t_1z^{32}\big)/\big(85
 z^{15}\big),\\
g_3 = \big(49t_1^6 + 2415t_1^4z^{16} + 3t_3^2z^{18} + 795t_1^2z^{32} - 35z^{48}\big)/
 \big(6z^{18}\big).
\end{gather*}
The determinant $\det(-T)$ regarded as a polynomial of $u_1$, $u_2$, $t_3$ (where $u_1=t_1/z^3$, $u_2=z^{10}$)
coincides with~$F_{B_6}$ by a weight preserving coordinate change up to a non-zero constant factor.

\subsection[Algebraic solutions related with the polynomial $F_{H_2}$ in \cite{Se1}]{Algebraic solutions related with the polynomial $\boldsymbol{F_{H_2}}$ in \cite{Se1}}\label{section6.5}

We recall the polynomial
\[
F_{H_2}=100x^3y^4+y^5+40x^4y^2z-10xy^3z+4x^5z^2-15x^2yz^2+z^3,
\]
which is a defining equation of a free divisor in $\mathbb{C}^3$ (cf.~\cite{Se1}).
There are two algebraic solutions which are related with the polynomial~$F_{H_2}$.

{\bf Icosahedral solution 29 of Boalch \cite{Bo1,Bo4} (LT18).} In this case,
\begin{gather*}
w(t_1)=\frac{1}{10},\qquad w(t_2)=\frac{1}{5},\qquad w(t_3)=1,\\
g_1 = -t_1\big(5t_1^6t_2^2 - 14t_2^5 - 2t_3\big)/2,\\
g_2 = \big(5t_1^{12} + 275t_1^6t_2^3 - 55t_2^6 + 33t_2t_3\big)/33, \\
g_3 = \big({-}100t_1^{18}t_2 + 2550t_1^{12}t_2^4 + 12750t_1^6t_2^7 + 595t_2^{10} + 9t_3^2\big)/18.
\end{gather*}
The determinant $\det(-T)$ regarded as a polynomial of $t_2$, $t_1^6$, $t_3$
coincides with $F_{H_2}$ by a weight preserving coordinate change up to a non-zero constant factor.

{\bf Icosahedral solution 30 of Boalch \cite{Bo1,Bo4} (LT19).}
In this case, $z$ is an algebraic function of $t_1$, $t_2$ defined by
\begin{gather*}
t_1^6 + t_2z^6 + z^9=0,\\
w(t_1)=\frac{3}{10},\qquad w(t_2)=\frac{3}{5},\qquad w(t_3)=1,\qquad w(z)=\frac{1}{5},\\
g_1=t_1\big({-}80t_2^2 + 910t_3z + 165t_2z^3 + 63z^6\big)/(910z), \\
g_2= \big(4t_2t_3 - 12t_2^2z^2 - 36t_2z^5 - 27z^8\big)/4,\\
g_3= \big(-560t_1^{18} + 595t_3^2z^{17} + 7140t_2^2z^{21} - 8160t_2z^{24} -
 15113z^{27}\big)/\big(1190z^{17}\big).
\end{gather*}
The determinant $\det(-T)$ regarded as a polynomial of $z$, $t_2$, $t_3$
coincides with $F_{H_2}$ by a weight preserving coordinate change up to a non-zero constant factor.

\subsection[Algebraic solution related with $E_{14}$-singularity]{Algebraic solution related with $\boldsymbol{E_{14}}$-singularity}\label{section6.6}

{\bf Octahedral solution 13 of Boalch \cite{Bo2} (LT30).} In this case,
\begin{gather*}
w(t_1)=\frac{1}{8},\qquad w(t_2)=\frac{3}{8},\qquad w(t_3)=1,\\
g_1=\big(5t_1^9 - 84t_1^6t_2 - 210t_1^3t_2^2 + 140t_2^3 + 9t_1t_3\big)/9,\\
g_2= \big(140t_1^{11} - 165t_1^8t_2 + 924t_1^5t_2^2 + 770t_1^2t_2^3 + 11t_2t_3\big)/11, \\
g_3= \big({-}95680t_1^{16} - 432320t_1^{13}t_2 + 780416t_1^{10}t_2^2 - 58240t_1^7t_2^3 + 1019200t_1^4t_2^4\\
\hphantom{g_3=}{} +
 203840t_1t_2^5 + 39t_3^2\big)/78.
\end{gather*}
The determinant $\det(-T)$ in this case coincides with
the polynomial
\begin{gather*}
F_{E_{14}} = -4x^6y^6 -\frac{ 20}{3}x^3y^7 - 3y^8 + 30x^7y^3z +
 51x^4y^4z + 24xy^5z - \frac{243}{4}x^8z^2\\
\hphantom{F_{E_{14}} =}{} - 108x^5yz^2 - 56x^2y^2z^2 - 8z^3
\end{gather*}
by a weight preserving coordinate change.
The polynomial $F_{E_{14}}$ is regarded as a 1-parameter deformation of the defining polynomial of $E_{14}$-singularity
in the sense of Arnol'd. In fact
\[
F_{E_{14}}|_{x=0}=-3y^8-8z^3.
\]
For topics related to $F_{E_{14}}$, see \cite{KMS, Se5}.

\subsection*{Acknowledgements}

Professor Yoshishige Haraoka taught the first author (M.K.)
 that integrable systems in three variables are useful to derive the Painlev\'e VI solutions.
 This is the starting point of our work.
 The authors would like to thank Professor Haraoka for his advice.
 After a preprint of this paper was written,
the authors received helpful comments including information on the papers
\cite{AL0,AL1,BoC,DH,Du2,Lo, Ma,Or} from Professors B.~Dubrovin, Y.~Konishi, C.~Hertling,
 P.~Boalch, A.~Arsie and P.~Lorenzoni, J.~Michel and H.~Terao.
 The authors express their sincere gratitude to these people.
 The authors thank anonymous referees for their useful comments and suggestions in order to improve the manuscript.
 This work was partially supported by JSPS KAKENHI Grant Numbers 25800082, 17K05335, 26400111, 17K05269.

\pdfbookmark[1]{References}{ref}
\LastPageEnding

\end{document}